
\input amstex
\documentstyle{amsppt}
\magnification=1200
\vsize19.5cm
\hsize13.5cm
\TagsOnRight
\pageno=1
\baselineskip=15.0pt
\parskip=3pt

\def\p{\partial}
\def\noo{\noindent}
\def\eps{\varepsilon}
\def\lam{\lambda}
\def\Om{\Omega}

\def\pom{{\p \Om}}

\def\R{\bold R}

\def\th{\theta}
\def\wtt{\tilde}

\def\Ga{\Gamma}

\def\det{\text{det}}

\def\ol{\overline}
\def\lan{\langle}
\def\ran{\rangle}
\def\D{\nabla}
\def\phi{\varphi}
\def\kk{k}

\def\M{\Cal M}

\def\F{\Cal F}
\def\E{\Cal E}

\nologo \NoRunningHeads

\topmatter

\title  {The Yamabe problem for higher order curvatures}\endtitle

\author{Weimin Sheng, Neil S. Trudinger, Xu-Jia Wang }\endauthor

\affil{   }\endaffil

\address Weimin Sheng, Department of Mathematics,
Zhejiang University, Hangzhou 310028, China      \endaddress
\email weimins\@zju.edu.cn \endemail

\address Neil S. Trudinger, Centre for Mathematics and Its Applications,
Australian National University,  Canberra, ACT 0200, Australia
\endaddress
\email neil.trudinger\@maths.anu.edu.au  \endemail

\address
Xu-Jia Wang, Centre for Mathematics and its Applications,
Australian National University, Canberra ACT 0200, Australia
\newline
Nankai Institute of Mathematics, Nankai University, Tianjin
300071, China
\endaddress
\email wang\@maths.anu.edu.au \endemail

\thanks
The first and third authors were supported by the Natural Science
Foundation of China. The second and third authors were supported
by the Australian Research Council
\endthanks

\abstract {Let $\M$ be a compact Riemannian manifold of dimension
$n$. The $\kk$-curvature, for $k=1,2, \cdots, n$, is defined as
the $k$-th elementary symmetric polynomial of the eigenvalues of
the Schouten tenser. The $\kk$-Yamabe problem is to prove the
existence of a conformal metric whose $\kk$-curvature is a
constant. When $k=1$, it reduces to the well-known Yamabe problem.
Under the assumption that the metric is admissible, the existence
of solutions to the $\kk$-Yamabe problem was recently proved by
Gursky and Viaclovsky for $k>\frac n2$. In this paper we prove the
existence of solutions for the remaining cases $2\le k \le \frac
n2$, assuming that the equation is variational.}
\endabstract


\endtopmatter

\document

\baselineskip=14.5pt
\parskip=3pt

\centerline{\bf 1. Introduction}

\vskip10pt

In recent years the Yamabe problem for the $\kk$-curvature of the
Schouten tensor, or simply the $\kk$-Yamabe problem, has been
extensively studied. Let $(\M, g_0)$ be a compact Riemannian
manifold of dimension $n$. Denote by $Riem$, $Ric$, and $R$ the
Riemannian curvature tensor, the Ricci tensor, and the scalar
curvature, respectively. Then one has the standard decomposition
$Riem =W +A\odot g_0$, where $W$ is the Weyl tensor, $A$ is the
Schouten tensor given in (1.2) below, and $\odot$ denotes the
Kulkarni-Nomizu product [B]. As the Weyl curvature tensor is
conformally invariant, the transformation of the Riemannian
curvature tensor under conformal changes of metrics is determined
by that of the Schouten tensor. Therefore it is of interest to
study curvature functions of the Schouten tensor under conformal
deformation. A fundamental problem is the $\kk$-Yamabe problem, to
prove the existence of a conformal metric
$g=g_v=v^{\frac{4}{n-4}}g_0$ whose $\kk$-curvature is equal to a
constant, that is
$$\sigma_k(\lam(A_g))=1,\tag 1.1$$
where $1\le k\le n$ is an integer, $\lam=(\lam_1, \cdots, \lam_n)$
are the eigenvalues of $A_g$ with respect to the metric $g$, and
$$A_g=\frac 1{n-2}(Ric_g-\frac {R_g} {2(n-1)}g) \tag 1.2$$
is the Schouten tensor. As usual we denote by
$$\sigma_k(\lam)
  =\sum_{i_1<\cdots<i_k}\lam_{i_1}\cdots\lam_{i_k}\tag 1.3$$
the $k$-th elementary symmetric polynomial. When $k=1$, we arrive
at the well known Yamabe problem.

When $k\ge 2$, the $\kk$-Yamabe problem was initiated by
Viaclovsky [V1] and also arose in the study of the Paneitz
operator [Br, CGY1]. Equation (1.1) is a fully nonlinear partial
differential equation. To work in the realm of elliptic operators,
one assumes that the eigenvalues $\lam(A_g)$ lie in the convex
cone $\Ga_k$ [CNS], where
$$\Ga_k=\{\lam\in \R^n\ |\
        \sigma_j(\lam)>0\ \text{for}\ j=1, \cdots, k\}. \tag 1.4$$
Under this assumption, the $\kk$-Yamabe problem has been solved in
the cases when $n=4$ and $k=2$ [CGY1, CGY2], or when the manifold
is locally conformally flat [LL1, GW2]. Very recently Gursky and
Viaclovsky [GV2] solved the problem for $k>\frac n2$, using the
positivity of the Ricci curvature in this case.

In this paper we employ a variational method to treat the problem
for the cases $2\le k\le \frac n2$. We prove that equation (1.1)
has a solution as long as it is variational, namely it is the
Euler equation of a functional, which includes the cases when
$k=2$ and when $\M$ is locally conformally flat. In Section 2 we
give a sufficient and necessary condition for equation (1.1) to be
variational.

When (1.1) is variational and $k\ne \frac n2$, its solutions
correspond to critical points of the functional
$$J(g)=\frac {n-2}{2(n-2k)}\int_\M \sigma_k(\lam(g))d\,vol_g
         -\frac {n-2}{2n}\int_\M d\,vol_g\tag 1.5$$
in the conformal class $[g_0]=\{g\ |\ g=v^{\frac{4}{n-2}}g_0,\
v>0\}$. When $k=\frac n2$, the first integral in (1.5) is a
constant and we need to replace it by (2.32) below. We will find a
min-max (Mountain Pass Lemma) solution, as in the case $k=1$. Note
that in the case $k>\frac n2$, the coefficients are negative and
the functional is negative.

The progressive resolution of the Yamabe problem $(k=1)$ by the
second author, Aubin and Schoen [Ya, Tr1, Au1, S1] was a milestone
in differential geometry. Roughly speaking, the overall proof
consists of two parts. The first one is to show that the Yamabe
problem is solvable if the Yamabe constant $Y_1$ satisfies the
condition
$$Y_1(\M)<Y_1(S^n),\tag 1.6$$
and the second one is to verify the condition (1.6) for manifolds
not conformally diffeomorphic to the unit sphere $S^n$ with
standard metric.  When $\M$ is locally conformally flat, different
proofs were found later [SY1, Ye].

For the $\kk$-Yamabe problem, $2\le k\le \frac n2$, our
variational approach basically comprises the same two steps.
Namely one first shows that (1.1) has a solution when the
$\kk$-Yamabe constant $Y_k$ satisfies
$$Y_k(\M)<Y_k(S^n),\tag 1.7$$
and then verify the condition (1.7) for manifolds not conformal to
the unit sphere $S^n$. But since equation (1.1) is fully
nonlinear, our treatment is technically different and more
complicated. For the first step, we cannot apply the variational
method directly, since we need to restrict the functional (1.5) to
a subset of the conformal class $[g_0]$, given by
$$[g_0]_k=\{g\ |\ g=v^{\frac{4}{n-2}}g_0,
                        v>0, \lam(A_g)\in\Ga_k\},\tag 1.8$$
and the set of functions $v$ with $v^{\frac{4}{n-2}}g_0\in
[g_0]_k$ may not be convex. Through the functional (1.5), we
introduce a descent gradient flow, establish appropriate a priori
estimates, and prove the convergence of solutions to the flow
under assumption (1.7). We need to choose a particular gradient
flow to obtain the a priori estimates, locally in time.

For the second step, it seems impossible to find an explicit
$k$-admissible test function. The function
$$v_\eps(x)=\big(\frac \eps{\eps^2+r^2}\big)^{(n-2)/2},\tag 1.9$$
which is the unique solution of (1.1) on the Euclidean space
$\R^n$ for all $1\le k\le n$, is $k$-admissible only when $r\le
C\eps^{1/2}$ on a general manifold, where $r$ is the geodesic
distance. Fortunately we found a simple way to deduce (1.7)
directly from (1.6).

This paper is arranged as follows. In Section 2,  we state the
main results, specifically in \S 2.1, while in \S 2.2 we outline
the proof. In \S 2.3 we collect some related results on the
$k$-Hessian equation. In \S 2.4 we give a sufficient and necessary
condition for a partial differential equation to be variational.
In Section 3 we study the regularity for the gradient flow of the
functional (1.5) for solutions with $\lam(A_{g_v})\in\Ga_k$, and
give counterexamples to interior regularity for solutions with
$\lam(A_{g_v})\in (-\Ga_k)$. In Section 4 we investigate the
asymptotic behavior of a descent gradient flow and prove the
convergence of the flow under condition (1.7). We then prove (1.7)
for manifolds not conformal to $S^n$ in Section 5. The final
Section 6 contains some remarks.

The authors are grateful to Kaiseng Chou for useful discussions.
This research was largely carried out in the winter of 2004-05
while the third author was at the Nankai Institute of Mathematics
in China under a Yangtze River Fellowship. The other authors are
also grateful for the Nankai Institute for hospitality when we
were all there together in November 2004.

\newpage

\vskip30pt

\centerline{\bf 2. The main results}

\vskip10pt

{\bf 2.1. The main results}. Let $(\M, g_0)$ be a Riemannian
manifold. If $g=v^{\frac {4}{n-2}}g_0$ is a solution of (1.1),
then the Schouten tensor is given by $A_g=\frac {2}{(n-2)v}V$, and
$v$ satisfies the equation
$$L[v] := v^{(1-k)\frac {n+2}{n-2}}\sigma_k(\lam(V))
       =v^{\frac{n+2}{n-2}}, \tag 2.1$$
where
$$V=-\D^2 v+\frac {n}{n-2}\frac {\D v\otimes\D v}{v}
    -\frac{1}{n-2}\frac{|\D v|^2}{v}g_0
    +\frac {n-2}{2} vA_{g_0}.  \tag 2.2$$

Equation (2.1) is a fully nonlinear equation of similar type to
the $k$-Hessian equations [CNS, CW2, I1, TW2]. For the operator
$L$ to be elliptic,  we need to restrict to metrics with
eigenvalues $\lam(A_g)\in \cup(\pm\Ga_k)$, which we will simply
denote as $g\in \pm\Ga_k$. Therefore equation (2.1) has two
elliptic branches, one is when the eigenvalues $\lam\in\Ga_k$ and
the other one is when $\lam\in(-\Ga_k)$. In this paper we will
mainly consider solutions with eigenvalues in $\Ga_k$. Accordingly
we say $v$ is {\it $k$-admissible} (that is $v$ is strictly
subharmonic with respect to $L$) if $g=v^{\frac{4}{n-2}}g_0 \in
\Gamma_k$. The set of all $k$-admissible functions will be denoted
by $\Phi_k=\Phi_k(\M, g_0)$. In this paper we will always assume,
unless otherwise indicated,  that $2\le k\le \frac n2$ and the
following two conditions hold,

\item {}{\bf (C1)} The set $\Phi_k(\M, g_0)\ne
\emptyset$;

\item {}{\bf (C2)} The operator $L$ is variational.

The condition (C1) ensures that the operator $L$ is elliptic, and
may be replaced by $Y_j(\M)>0$ for $j=1, \cdots, k$, as in the
case when $k=2$ and $n=4$ [CGY1, GV1]. Note that in condition
(C1), we do not assume directly that the metric $g_0\in \Ga_k$,
rather we assume that there exists a positive function $v$ such
that $v^{\frac {4}{n-2}}g_0\in\Ga_k$. Conditions (C1) (C2) are
automatically satisfied when $k=1$.

As for the Yamabe problem, we introduce the $\kk$-Yamabe constant
for $2\le k\le  \frac n2 $,
$$Y_k(\M)=\inf
   \{\F_k(g)\ |\ g\in [g_0]_k, \text{Vol}(\M_g)=1\}, \tag 2.3$$
where $[g_0]_k$ is defined in (1.8),  and
$$\align
\F_k(g)  & =\int_\M \sigma_k(\lam(A_g))d\, vol_g\\
 &=\int_{\M}   v^{\frac {2n}{n-2}-k\frac {n+2}{n-2}}
                       \sigma_k(\lam(V))\,d\, vol_{g_0}.\tag 2.4\\
  \endalign $$
Note that we have ignored a coefficient $(\frac 2{n-2})^k$ in the
second equality. The main result of the paper is the following.

\proclaim{Theorem 2.1} Assume $2\le k\le\frac n2$ and the
conditions (C1) (C2) hold. Then the $\kk$-Yamabe problem (1.1) is
solvable.
\endproclaim

As indicated in the introduction, the proof of Theorem 2.1 is
divided into two parts. The first part is the following lemma.

\proclaim{Lemma 2.1} If the critical inequality (1.7) holds, then
the $k$-Yamabe problem (1.1) is solvable.
\endproclaim

The second part provides the condition for (1.7).

\proclaim{Lemma 2.2}
 The critical inequality (1.7) holds for any compact manifold
 which is not conformal to the unit sphere $S^n$.
\endproclaim

When $k=\frac n2$, we prove that $\F_{n/2}(g)\equiv Y_{n/2}(\M)$,
that is it is a constant for any $g\in [g_0]_k$ (Lemma 4.8). Hence
(1.7) implies that $\F_{n/2}(g)<Y_{n/2}(S^n)$ provided $\M$ is not
conformal to the unit sphere.

\vskip10pt

{\bf 2.2. Strategy of the proof}. A solution of the $\kk$-Yamabe
problem is a min-max type critical point of the corresponding
functional. As we need to restrict ourselves to $k$-admissible
functions, we cannot directly use variational theory (such as the
Ekeland variational principle). Instead we study a descent
gradient flow of the functional and investigate its convergence.
We need to choose a special gradient flow (similar to [CW2]) for
which the necessary a priori estimates can be established.

As with the original Yamabe paper [Ya], we first study the
approximating problems
$$L(v) =v^p, \tag 2.5$$
where $1<p\le\frac {n+2}{n-2}$. Equation (2.5) is the Euler
equation of the functional
$$J_p(v)=J_p(v; \M)=\frac{n-2}{2n-4k}\int_{(\M, g_0)}
       v^{\frac {2n}{n-2}-k\frac {n+2}{n-2}}\sigma_k(\lam(V))
              -\frac {1}{p+1}\int_{(\M, g_0)} v^{p+1}. \tag 2.6$$

Let $\phi_1=\eps$ and $\phi_2=\eps^{-1}$, where $\eps>0$ is a
small constant. Then $J_p(\phi_1)\to 0$ (when $k<\frac n2$) and
$J_p(\phi_2)\to-\infty$ as $\eps\to 0$. Let $P$ denote the set of
paths in $\Phi_k$ connecting $\phi_1$ and $\phi_2$, namely
$$P=\{\gamma\in C([0, 1], \Phi_k)\ |\
      \gamma(0)= \phi_1, \gamma(1)=\phi_2\}.\tag 2.7$$
Obviously $\Phi_k\ne\emptyset$. Denote
$$c_p[\M] =\inf_{\gamma\in P}\sup_{s\in [0, 1]}J_p(\gamma(s);\M).
                                                \tag 2.8$$
Then (1.7) is equivalent to
$$c_p[\M]<c_p[S^n] \tag 2.9$$
with $p=\frac {n+2}{n-2}$. We will prove that $J_p$ has a min-max
critical point $v_p$ with $J_p(v_p)=c_p[\M]$,  in the sub-critical
case $p<\frac {n+2}{n-2}$. By a blow-up argument, we prove
furthermore that $v_p$ converges to a solution of (2.1) under the
assumption (2.9).

The descent gradient flow will be chosen so that appropriate a
priori estimates can be established. To simplify the computations,
we will also use the conformal transformations $g=u^{-2}g_0$ or
$g=e^{-2w}g_0$. That is
$$u=e^w=v^{-\frac {2}{n-2}}. \tag 2.10$$
We say $u$ or $w$ is $k$-admissible if $v$ is, and also denote $u,
w\in \Phi_k$ if $v\in \Phi_k$.

Our gradient flow is given by
$$F[w]-w_t=\mu(f(x, w)), \tag 2.11$$
where
$$F[w] : =\mu(\sigma_k(\lam(A_g )))    \tag 2.12$$
and $g=e^{-2w}g_0$. When $f(x, w)=e^{-2kw}$, a stationary solution
of (2.11) is a solution to the $\kk$-Yamabe problem. The function
$\mu$ is monotone increasing and satisfies
$$\lim_{t\to 0^+} \mu(t)=-\infty .\tag 2.13$$
Condition (2.13) ensures the solution is $k$-admissible at any
time $t$. For if $u(\cdot, t)$ is a smooth solution, then (2.13)
implies $\sigma_k(\lam)>0$ at any time $t>0$.  A natural candidate
for the choice of $\mu$ is the logarithm function $\mu(t)=\log t$
[Ch1, W1, TW4]. However for the flow (2.11), we need to choose a
different $\mu$ to ensure appropriate a priori estimates.

In the case $k=\frac n2$, $\F_{n/2}(g)$ is a constant less than
$Y_{n/2}(S^n)$. Hence by the Liouville theorem in [LL1], it is
easy to prove that the set of solutions of (2.1) is compact. Hence
the existence of solutions can be obtained by a degree argument.
When $2\le k<\frac n2$, by the Liouville theorem in [LL2], one can
also prove the set of solutions of (2.5) is compact when
$p<\frac{n+2}{n-2}$. But to use the condition (1.7) in the blow-up
argument, we need a solution $v_p$ of (2.5) satisfying
$J_p(v_p)=c_p$. This is the reason for us to employ the gradient
flow.

\vskip10pt

For the verification of (1.7), one cannot mimic the argument in
the case $k=1$, as the test function (1.9) is in general not
$k$-admissible in a geodesic ball $B_{\rho_0}$. Instead we let
$v_1$ be the solution to the Yamabe problem ($k=1$), and $v$ be
the $k$-admissible solution of the equation
$$\sigma_k(\lam(V))=v_1^{k\frac{n+2}{n-2}} .$$
We will verify (1.7) by using the solution $v$ as the test
function.

\vskip10pt

The idea of using a gradient flow was inspired by [CW2], where a
similar problem for the $k$-Hessian problem (see (2.22) below) was
studied. However technically the argument in this paper is
different. For the $\kk$-Yamabe problem, the corresponding Sobolev
type inequality was not available, and the a priori estimates only
allow us to get a local (in time) solution. The argument in this
paper is also self-contained, except we will use the Liouville
theorem in [LL1, LL2], proved by the moving plane method; see also
[CGY3].

\vskip20pt

{\bf 2.3. The $k$-Hessian equation}. Equations (2.1) is closely
related to the $k$-Hessian equation
$$\sigma_k(\lam(D^2 v))=f(x)\ \ \ x\in\Om, \tag 2.15$$
where $1\le k\le n$, $\lam=(\lam_1, \cdots, \lam_n)$ denote the
eigenvalues of the Hessian matrix $(D^2 v)$, $\Om$ is a bounded
domain in the Euclidean $n$-space $\R^n$. For later applications
we collect here some elementary properties of the polynomial
$\sigma_k$, and give a very brief summary of related results for
the equation (2.15).

We write $\sigma_0(\lam)=1$, $\sigma_k(\lam)=0$ for $k>n$, and
denote $\sigma_{k; i}(\lam)=\sigma_k(\lam)_{\big{|}\lam_i=0}$.

\proclaim {Lemma 2.3} Let $\lam\in\Ga_k$ with $\lam_1\ge
\cdots\ge\lam_n$. Then
$$\align
&\lam_k\ge 0\tag i\\
&\sigma_k(\lam)
=\sigma_{k;i}(\lam)+\lam_i\sigma_{k-1;i}(\lam),\tag ii\\
&\Sigma_{i=1}^n \sigma_{k-1;i}(\lam)
                    =(n-k+1)\sigma_{k-1}(\lam),\tag iii\\
&\sigma_{k-1; n}(\lam)\ge\cdots\ge\sigma_{k-1; 1}(\lam)>0,\tag iv\\
&\sigma_{k-1;k}(\lam)\ge
           C_{n, k}\Sigma_{i=1}^n \sigma_{k-1;i}(\lam),\tag v\\
&\sigma_{k-1}(\lam)\ge \frac{k}{n-k+1}(^n_k)^{1/k}
[\sigma_k(\lam)]^{(k-1)/k}. \tag vi\\
\endalign $$
Moreover, the function $[\sigma_k]^{1/k}$ is concave on $\Ga_k$.
\endproclaim

We just listed a few basic formulae, there are many other useful
ones, see for example [CNS, CW2, LT]. For our investigation of
equation (2.1) and its parabolic counterpart, Lemma 2.1 will be
sufficient. These formulas can be extended to $\sigma_k(\lam(r))$,
regarded as functions of $n\times n$ symmetric matrices $r$. In
particular $[\sigma_k(\lam(r))]^{1/k}$ is concave in $r$ [CNS].

We say a function $v\in C^2(\Om)$ is $k$-admissible (relative to
equation (2.15)) if the eigenvalues $\lam(D^2 v)\in \Ga_k$.
Equation (2.15) is elliptic if $v$ is $k$-admissible. The
existence of $k$-admissible solutions to the Dirichlet problem for
(2.15) was proved by Caffarelli-Nirenberg-Spruck [CNS], see also
Ivochkina [I].

Relevant to the $\kk$-Yamabe problem is the variational property
of the $k$-Hessian equation (2.15), investigated in [CW2, TW4,
W1]. It is well known that the $k$-Hessian equation is the Euler
equation of the functional
$$I_k(v)=\frac 1{k+1}\int_\Om (-v) \sigma_k(\lam(D^2 v)). \tag 2.16$$
The Sobolev-Poincar\'e type inequality, for $k$-admissible
functions vanishing on the boundary,
$$I_l^{1/(l+1)}(v)\le CI_k^{1/(k+1)}(v), \tag 2.17$$
was established in [W1] for the case $l=0$ and $k\ge 1$, and in
[TW4] for the case $k>l\ge 1$, where $0\le l\le k\le n$,
$$I_0(v)=\big[\int_\Om |v|^{k^*} dx\big]^{1/k^*}, \tag 2.18$$
and
$$\align
k^* & =n(k+1)/(n-2k)\ \ \ \text{if}\ \ k<n/2,\tag 2.19\\
k^* & <\infty \ \ \ \text{if}\ \ k=n/2,\\
k^* & =\infty\ \ \ \text{if}\ \ k>n/2.\\
\endalign$$
The best constant in the inequality (2.17) is attained by
$$v(x)= (1+|x|^2)^{(2k-n)/2k}\tag 2.20$$
when $l=0$, $k<\frac n2$, and $\Om=\R^n$; and by the unique solution
of
$$\frac {\sigma_k}{\sigma_l}(\lam(D^2 v))=1
                         \ \ \text{in}\ \ \Om\tag 2.21$$
for $1\le l<k\le n$.

>From the inequalities (2.17), it was proved in [CW2] that the
Dirichlet problem
$$\align
\sigma_k(\lam(D^2 v)) & =|v|^p+f(v)
      \ \ \ \text{in}\ \ \Om,\tag 2.22\\
      v & =0\ \ \ \text{on}\ \ \pom,\\
      \endalign $$
admits a nonzero $k$-admissible solution, where $1\le k\le \frac
n2$, $1<p< k^*-1$, $f$ is a lower order term of $|v|^p$. The
existence result was proved for the problem with more general
right hand side.

\vskip5pt

\noo{\bf Remark 2.1}. The 1996 preprint [CW1] also contains the
existence of solutions to the Dirichlet problem (2.22) in the
critical growth case $p=k^*-1$. The result was obtained by a
blow-up argument and the symmetrization of functions, see Theorem
9.1 in [CW1].

\vskip20pt

\noo {\bf 2.4. A necessary and sufficient condition for an
equation to be variational}.

The following proposition was communicated to the authors by
Kaiseng Chou several years ago.

\noo\proclaim{Proposition 2.1} Let $\M$ be a compact manifold
without boundary, $v\in C^4(\M)$. An operator $F[v]=F[\D ^2v, \D
v, v, x]$ is variational if and only if its linearized operator is
self-adjoint. The functional is given by
$$I[v]=\int G[v],\tag 2.23$$
except when $F$ is homogeneous of degree -1,  where
$$G[v]=\int_0^1 v F[\lam v].\tag 2.24 $$
\endproclaim

This proposition can be found in [O].  We give a proof of the
$``$if$"$ part, as we need some related formulae.

{\it Proof.} The linearized operator of $F[v]$ is given by
$$L(\phi)  =F^{ij}\phi_{ij}+F_{p_j}\phi_j +F_v \phi.\tag 2.25$$
We have
$$\align
\int_{\M} v\, L(\phi)
 &=\int_{\M} [v\D _i(F^{ij}\D _j\phi)+ v\phi F_v] -A \\
 &=\int_{M} [-v_i\phi_jF^{ij}+v\phi F_v] -A\\
 &=\int_{\M} \phi [F^{ij} v_{ij}+F_{p_i}v_i+F_vv] -A+B\\
 &=\int_{\M} \phi\,L(v) -A+B,\\
\endalign $$
where
$$\align
A&=\int_{\M} v\phi_j(\D _iF^{ij}-F_{p_j}),\\
B&=\int_\M v_i\phi(\D _jF^{ij}-F_{p_i}),\\
-A+B& = -\int_\M
     (\frac{\phi}{v})_i v^2(\D _jF^{ij}-F_{p_i})\\
     & = \int_\M
     \frac{\phi}{v}\D _i[ v^2(\D _jF^{ij}-F_{p_i})]\\
\endalign $$
Hence $L$ is self-adjoint if and only if
$$\sum_{i, j=1}^n \D _i[ v^2(\D _jF^{ij}-F_{p_i})]=0. \tag 2.26$$

If $L$ is self-adjoint,
$$\align
\lan I'[v], \phi\ran
 & =\int_\M\phi\int_0^1 F[\lam v]
 +\int_0^1\int_\M \lam v[F^{ij}[\lam v]\phi_{ij}
          +F_{p_i}\phi_i+F_v\phi]\\
 & =\int_\M\phi\int_0^1 F[\lam v]
  +\int_0^1\int_\M \lam \phi [F^{ij}[\lam v]v_{ij}
          +F_{p_i}v_i+F_vv]\\
 & =\int_\M\phi\int_0^1 F[\lam v]
 +\int_\M\phi \int_0^1\lam\frac {d}{d\lam} F[\lam v]d\lam\\
 &= \int_\M \phi F[v]. \\
 \endalign $$
Hence $F$ is the Euler equation of the functional $I$. $\square$

Conversely if the operator $F$ is the Euler operator of the
functional $I$, from the above argument we must have $-A+B=0$,
namely (2.26) holds. In other words, $F$ is the Euler operator of
$I$ if and only if (2.26) holds. Observe that if
$$\sum_i \D _iF^{ij} =F_{p_j}\ \ \ \forall\ j,\tag 2.27$$
then (2.26) holds.

>From Proposition (2.5) we can recover the results on the
variational structure of (2.1) in [V1]. First, if locally $\M$ is
Euclidean, one verifies directly that (2.26) holds, as it is a
pointwise condition. The locally conformally flat case is
equivalent to the Euclidean case by a conformal deformation to the
Euclidean metric. Finally if $k=2$, we note that to verify (2.26)
for arbitrary $v$ with a fixed background metric $g_0$ is
equivalent to verify it for $v\equiv 1$ with respect to an
arbitrary conformal metric $g={\hat v}^{\frac{4}{n-2}}g_0$.
However when $v\equiv 1$, condition (2.26) becomes $\sum_{i,j=1}^n
\D_i\D_j F^{ij}=0$, where $F^{ij}=\frac {\p}{\p
r_{ij}}\sigma_k(\lam(r))$ at $r=A_g$. But we have
$$\nabla_i F^{ij}=\frac{1}{2(n-2)}(R,_j-2R_{ij,i})=0 \tag 2.28$$
by the second Bianchi identity.

When $k\ge 3$, it is easy to find metrics for which (2.26) does
not hold (at $v\equiv 1$). So equation (2.1) need not be
variational. As an example, let $k=n=3$, and in a local coordinate
system, let the metric $g=\{g_{ij}\}$ be given by
$$g_{11}=1,\ \ g_{22}=1+x,\ \ g_{33}=1+y^2+z^2,
 \ \ g_{ij}=\delta_{ij}\ \ \text{for}\ \ i\neq j.\tag 2.29$$
><DEFANGED.171172 Then $\nabla_i\nabla_j F^{ij}\neq 0$.

By (2.23) we also see that (2.6) is the functional of (2.5). When
$k=\frac n2$, the integral (2.24) may not exist. We may consider
$v$ as a composite function $v=\phi(w)$ and write equation (2.5)
in the form
$$F[w]=: \phi'(w)L(\phi(w))=\phi^p(w)\phi'(w).\tag 2.30$$
If the operator $L$ in (2.1) satisfies (2.26) with respect to $v$,
the operator $F$ in (2.30) satisfies (2.26) with respect to $w$.
Hence the corresponding functional is given by
$$\E_{n/2}(w)=\int_{(\M, g_0)}\int_0^1 w F[tw]. \tag 2.31$$
In particular if $v=e^{-\frac {n-2}{2}w}$, then we obtain the
functional in [BV],
$$\E_{n/2}(w) =-\int_{(\M, g_0)}\int_0^1
   w \sigma_{n/2}(\lam(A_{g_t})),\tag 2.32$$
where $g_t=e^{-2tw}g_0$.

\newpage

\vskip30pt

\centerline{\bf 3. The a priori estimates}

\vskip10pt

In this section we study the regularity of $k$-admissible
solutions ($2\le k\le n$) to equation (2.1) and its parabolic
counterpart (2.11). The global a priori estimates for the elliptic
equation (2.1) (for solutions with eigenvalues in $\Ga_k$) have
already been established by J. Viaclovsky [V2], with interior
estimates by P. Guan and G. Wang [GW1]. We will provide a simpler
proof for the elliptic equation (2.1), and extend the estimates to
the parabolic equation (2.11) on general manifolds, which is
necessary for our proof of Theorem 2.1. Previously the estimates
for the parabolic equation were proved on locally conformally flat
manifolds [Ye, GW2, GW3]. Regularity has also been studied in many
other papers [CGY1,LL1].

We will also present an example showing that the interior a priori
estimates do not hold for solutions with eigenvalues in the
negative cone $-\Ga_k$.

\vskip10pt

{\bf 3.1. A priori estimates for equation (2.1)}. For the regularity
of (2.1), we will use the conformal changes $g=u^{-2}g_0$. For
function $u$, equation (2.1) becomes
$$\sigma_k(\lam(U))=u^{-k}, \tag 3.1$$
where
$$U=\D^2 u-\frac {|\D u|^2}{2u}g_0+uA_{g_0}. $$

\proclaim{Lemma 3.1} [GW1] Let $u\in C^3$ be a $k$-admissible
positive solution of (3.1) in a geodesic ball $B_{r}(0)\subset\M$.
Suppose $A_{g_0}=(a_{ij})\in C^1(B_{r}(0))$. Then we have
$$\frac {|\D u|}{u}(0)\le C,\tag 3.2$$
where $C$ depends only on $n, k$, $r$, $\inf u$, and
$\|A_{g_0}\|_{C^1}$, $\D $ denotes the covariant derivative with
respect to the initial metric $g_0$.
\endproclaim

\noo{\it Proof}. Let $\mu$ be a smooth, monotone increasing
function. Write equation (3.1) in the form
$$F[u] =\mu[f(x, u)],  \tag 3.3$$
where $F[u]=\mu[\sigma_k(\lam(U))]$.  We will prove (3.2) for more
general function $f$. Moreover the constant $C$ is independent of
$\sup_{B_r} u$ if $f=\kappa u^{-p}$ for some constant $p>0$ and
smooth, positive function $\kappa$.

Let $z=|\D u|^2\phi^2(u)\rho^2$, where $\phi(u)=\frac 1u$, and
$\rho(x)=(1-\frac {|x|^2}{r^2})^+$ is a cut-off function, $|x|$
denotes the geodesic distance from $0$. Suppose $z$ attains
maximum at $x_0\in B_1(0)$, and $|\D u(x_0)|=u_1(x_0)$. Then at
$x_0$, in an orthonormal frame,
$$\align
\frac 12 (\log z)_i  &= \frac {u_{1i}}{u_1}
       +\frac {\phi'}{\phi} u_i+\frac {\rho_i}{\rho}=0, \tag 3.4 \\
\frac 12 (\log z)_{ij} &= \frac {u_{1ij}}{u_1}
       +\sum_{\alpha>1}\frac {u_{\alpha i}u_{\alpha j}}{u_1^2}
  -\frac {u_{1i}u_{1j}}{u_1^2}+\frac {\phi'}{\phi} u_{ij}
       +(\frac {\phi''}{\phi}-\frac {{\phi'}^2}{\phi^2})u_iu_j
       +(\frac {\rho_{ij}}{\rho}-\frac {\rho_i\rho_j}{\rho^2}).
                                             \tag 3.5\\
                     \endalign $$
Differentiating equation (3.3) gives
$$F^{ij}[u_{ij1}+(\frac {u_1^3}{2u^2}
       -\frac {u_1u_{11}}{u}) \delta_{ij}]=\Delta ,\tag 3.6$$
where for a matrix $r=(r_{ij})$, $F^{ij}(r)=\frac {\p}{\p r_{ij}}
\mu[\sigma_k(\lam(r))] =\mu'\frac{\p}{\p r_{ij}}\sigma_k(\lam(r))$,
$$\Delta=\D_1 \mu(f)
                   -F^{ij}\D_1 (a_{ij}u).$$
By (3.4)-(3.6) we have, at $x_0$,
$$\align
0&\ge \frac 12 F^{ij}(\log z)_{ij}
 = \frac{1}{u_1}[\frac {u_1u_{11}}{u}-\frac{u_1^3}{2u^2}]\Cal F
   +\sum _{\alpha>1}F^{ij}\frac {u_{\alpha i}u_{\alpha j}}{u_1^2}
  -F^{ij}(\frac {\phi'}{\phi}u_i+\frac {\rho_i}{\rho})
          (\frac {\phi'}{\phi}u_j+\frac {\rho_j}{\rho})\\
 & +\frac {\phi'}{\phi}F^{ij}(u_{ij}-\frac{|\D u|^2}{2u}\delta_{ij})
   + \frac {u_1^2}{2u}\frac {\phi'}{\phi} \Cal F
 +(\frac {\phi''}{\phi}-\frac {{\phi'}^2}{\phi^2})F^{11}u_1^2
  +F^{ij}(\frac {\rho_{ij}}{\rho}-\frac {\rho_i\rho_j}{\rho^2})
   +\frac{\Delta}{u_1}+\Delta', \\
\endalign $$
where $\Cal F=\sum F^{ii}$, $\Delta'$ arises in the exchange of
derivatives, with $|\Delta'|\le C\Cal F$. Note that
$$F^{ij}(u_{ij}-\frac{|\D u|^2}{2u}\delta_{ij})
 = k\mu'\sigma_k(\lam)-uF^{ij}a_{ij}\ge -C_au\Cal F, $$
where $C_a=0$ if $A_{g_0}=(a_{ij})=0$. By (3.4) and since
$\phi(u)=\frac 1u$,
$$\align
[\frac {u_{11}}{u}-\frac{u_1^2}{2u^2}]
       +\frac {u_1^2}{2u}\frac {\phi'}{\phi}
 & =-\frac {u_1\rho_1}{u\rho},\\
 -F^{ij}(\frac {\phi'}{\phi}u_i+\frac {\rho_i}{\rho})
          (\frac {\phi'}{\phi}u_j+\frac {\rho_j}{\rho})
  & +(\frac {\phi''}{\phi}-\frac {{\phi'}^2}{\phi^2})F^{11}u_1^2
   =-F^{ij}(\frac {2\phi' u_i\rho_j}{\phi \rho}+\frac
{\rho_i\rho_j}{\rho^2} )\\
\endalign $$
Hence we obtain
$$0\ge
 \sum_{\alpha>1} F^{ij}\frac {u_{\alpha i}u_{\alpha j}}{u_1^2}
 -C(\frac {1}{r^2\rho^2} +\frac{u_1}{u}\frac{1}{r\rho}+C_a)\Cal F
                    +\frac{\Delta}{u_1}+\Delta'. \tag 3.7$$
Denote $b=\frac {|\D u|^2}{2u}(x_0)$. We claim
$$\sum_{\alpha>1} F^{ij}u_{\alpha i}u_{\alpha j}
                    \ge Cb^2\Cal F-C'u^2\Cal F\tag 3.8$$
for some positive constant $C, C'$ ($C'=0$ if $a_{ij}=0$). Note
that by Lemma 2.3 (iii) (vi),  $\Cal F\ge C_{n, k}\mu'
\sigma_k^{(k-1)/k}$. From (3.8) we have
$$\frac{|\D u|}{u}\rho\le \frac {C_1}{r}+C_2\tag 3.9$$
at $x_0$, where $C_1$ is independent of $f$ and $C_2$ independent
of $r$. Hence $z(0)\le z(x_0)\le C$, namely (3.2) holds.

Denote $\wtt u_{ij}=u_{ij}+u a_{ij}$. For any two unit vectors
$\xi, \eta$, we denote formally $\wtt u_{\xi\eta}=\sum
\xi_i\eta_j\wtt u_{ij}$. Then to prove (3.8) it suffices to prove
$$A=:\sum_{\alpha>1} F^{ij}{\wtt u}_{\alpha i}{\wtt u}_{\alpha j}
                        \ge Cb^2\Cal F.\tag 3.10$$
By a rotation of the coordinates we suppose $\{\wtt u_{ij}\}$ is
diagonal at $x_0$. Then
$$\lam_1=\wtt u_{11}-b, \cdots, \lam_n=\wtt u_{nn}-b$$
are the eigenvalues of the matrix $\{\wtt u_{ij}-\frac
{|Du|^2}{2u}\delta_{ij}\}$. Suppose $\lam_1\ge \cdots\ge\lam_n$. At
$x_0$ we have $|Du(x_0)| =u_\xi(x_0)$ for some unit vector $\xi$. In
the new coordinates we have
$$A=\sum_i(F^{ii}{\wtt u}_{ii}^2-F^{ii}{\wtt u}_{\xi i}^2). $$

If there exists a small $\delta_0>0$ such that $\lan e_i,
\xi\ran<1-\delta_0$ for all unit axial vectors $e_i$, then $A\ge
\delta_0 F^{ii}\wtt u_{ii}^2$. Since $\lam=(\lam_1, \cdots,
\lam_n)\in\Ga_k^+$, we have $\lam_k>0$ and so $\wtt u_{kk}>b$.
Hence by Lemma 2.3(v), $A\ge \delta_0 b^2 F^{kk}\ge
\delta_1b^2\Cal F$. We obtain (3.8).

If there is $i^*$ such that $\lan e_{i^*}, \xi\ran\ge 1-\delta_0$,
then we have $A\ge\frac 12 \sum_{i\ne i^*} F^{ii}\wtt u_{ii}^2$.
If there exists $j\ge k$, $j\ne i^*$ such that $\wtt u_{jj}\ge
\alpha b$ for some $\alpha>0$, then by Lemma 2.3(iv)(v), $A\ge
\frac 12 F^{jj}(\alpha b)^2\ge\delta_2 b^2\Cal F $ and the claim
holds. Otherwise we have $i^*=k$ since $\wtt u_{kk}=\lam_k+b\ge
b$.

Case 1: $k\le n-2$. Observing that $\frac {\p}{\p\lam_1}\cdots
\frac{\p}{\p\lam_{k-1}}\sigma_k(\lam)=\lam_k+\cdots+\lam_n\ge 0$, we
have $\lam_k\ge -(\lam_{k+1}+\cdots+\lam_n)$. Since $\wtt u_{jj}\le
\alpha b$ for $j\ge k+1$, we have $\lam_j\le -(1-\alpha)b$. Hence
$\lam_k \ge (n-k)(1-\alpha) b\ge 2(1-\alpha)b$.

On the other hand, by (3.4), we may suppose at $x_0$,
$|\frac{\rho_\xi}\rho|\le \alpha \frac {u_\xi}{u}$, for otherwise
we have the required estimate (3.2). Hence $\wtt u_{\xi\xi}\le
(2+\alpha)b$ for a different small $\alpha>0$. By the relation
$\wtt u_{\xi\xi}=\sum_i \xi^2_i\wtt u_{ii}\ge \sum_{i\le
k}\xi^2_i\wtt u_{ii}-n\alpha b$ where $\xi=(\xi_1, \cdots,
\xi_n)$, we have $\wtt u_{kk}\le (1+\alpha)\wtt u_{\xi\xi}\le
(2+2\alpha)b$. Hence $\lam_k=\wtt u_{kk}-b\le (1+2\alpha)b$. We
reach a contradiction when $\alpha$ is sufficiently small.

Case 2: $k=n-1$.  We have
$$\frac {\p\sigma_k}{\p\lam_{k-1}} \lam_{k-1} = \sigma_k(\lam)
            -\frac {\lam_1\cdots\lam_n}{\lam_{k-1}}
               \ge -\frac {\lam_1\cdots\lam_n}{\lam_{k-1}}.$$
Since $\lam_n=\wtt u_{nn}-b\le -(1-\alpha) b$ and by $\frac
{\p}{\p\lam_1}\cdots\frac{\p}{\p\lam_{n-2}}\sigma_k(\lam)
=\lam_{n-1}+\lam_n\ge 0$, we have $\lam_{n-1}\ge (1-\alpha)b$ and
so $\lam_i\ge (1-\alpha)b$ for any $1\le i\le n-1$.  Hence $\frac
{\p\sigma_k}{\p\lam_{k-1}} \lam_{k-1}\ge (1-\alpha)^2
b^2\lam_1\cdots\lam_{n-3}$. Note that
$\mu'\frac{\p\sigma_k}{\p\lam_i}(\lam)=F^{ii}$. It follows that
$$\align
 A&\ge \frac 12\mu'
       \frac {\p\sigma_k}{\p\lam_{k-1}}\wtt u^2_{k-1\,k-1}
  \ge \frac 12\mu'
       \frac {\p\sigma_k}{\p\lam_{k-1}}\lam_{k-1}^2\tag 3.11\\
  &\ge \frac 12\mu' (1-\alpha)^2 b^2 \lam_1\cdots\lam_{n-2}
  \ge Cb^2\Cal F.\\
  \endalign$$

Case 3: $k=n$. As in Case 1, we assume that $|\frac{\rho_\xi}\rho|
\le \alpha \frac {u_\xi}{u}$. Then by (3.4), $\wtt u_{\xi\xi}\ge
(2-\alpha)b$. Note that when $k=n$, $\lam_i>0$ for all $i$. Hence
$\wtt u_{ii}=\lam_i+b>b$. Recall that when $k=n$, we have $i^*=n$.
It follows that $\wtt u_{nn}\ge (2-\alpha)b$ for a different small
$\alpha>0$. Hence $\lam_n\ge (1-\alpha)b$ and
$$A \ge \frac 12 \Sigma_{i\ne i^*} F^{ii}\wtt u_{ii}^2
    \ge \frac 12 F^{ii}\lam_i^2
    =   \frac 12 \lam_i\lam_n F^{nn}\ge Cb^2\Cal F.\tag 3.12$$
This completes the proof. $\square$

We remark that in our proof of (3.10), we didn't use the equation
(3.3). Hence we can also use (3.10) for the corresponding equation
(3.21) below. We also note that the gradient estimate is
independent of the choice of $\mu$. From Lemma 3.1, we obtain the
following Liouville theorem.

\proclaim{Corollary 3.1} Let $u\in C^3$ be an entire
$k$-admissible positive solution of
$$\sigma_k \big(\lam(\D^2 u -\frac {|\D u|^2}{2u}I)\big)
                                        =0.\tag 3.13$$
Then $u\equiv constant$.
\endproclaim

\noo{\it Proof}. For equation (3.13), the constant $C_2$ in (3.9)
vanishes. Letting $r\to\infty$, by (3.9), we see that either
$\frac {|\D u|}{u}\equiv 0$, or $\Cal F=0$. In the former case,
$u$ is a constant. In the latter case, $u$ satisfies
$\sigma_{k-1}(\lam)=0$ and so it is also a constant by induction.
$\square$

By approximation as in [MTW,W2] one can show that Corollary 3.1
holds for continuous positive viscosity solutions. The proof of
the interior gradient estimate (3.2) can be simplified if one
allows the estimate to depend on both $\inf_{B(0, r)} u$ and
$\sup_{B(0, r)} u$. Indeed, let $\phi(u)=\frac {1}{u-\delta}$ in
the auxiliary function $z$, where $\delta=\frac 12\inf_{B(0, r)}
u$. Then one obtains the extra good term $\frac {\delta
u_1^2}{(u-\delta)u^2} \Cal F$ on the right hand side of (3.7). The
proof after (3.8) is not needed.

For the $\kk$-Yamabe problem, $f(u)=u^{-k}$. The constant $C$ in
(3.2) is independent of $\sup u$. Therefore we have the Harnack
inequality [GW1].

\proclaim{Corollary 3.2} Let $u\in C^3$ be a positive solution of
(3.1). If $\inf u\ge C_0>0$, then $\sup u\le C_1$.
\endproclaim

Next we prove the second order derivative estimate.

\proclaim{Lemma 3.2} Let $u\in C^4$ be a $k$-admissible positive
solution of (3.1) in a geodesic ball $B_r(0)\subset\M$. Suppose
$A\in C^2(B_r(0))$. Then we have
$$|\D^2u|(0)\le C,\tag 3.14$$
where $C$ depends only on $n, k$, $r$, $\inf u$, $\sup u$, and
$\|A_{g_0}\|_{C^2}$.
\endproclaim

\noo{\it Proof}. Again we will consider the more general equation
(3.3). Choose $\mu(t)=t^{1/k}$ such that equation (3.3) is concave
in $U_{ij}$. Differentiating (3.3) we get
$$  F^{ij} U_{ij, kk}   =-\frac {\p^2 \mu(\sigma_k(\lam(U)))}
  {\p U_{ij}\p U_{rs}}U_{ij,k}U_{rs,k}+\D_k^2 \mu(f)
    \ge \D_k^2 \mu(f) ,\tag 3.15 $$
where $U_{ij, k}=\D_k U_{ij}$. As above denote $\wtt
u_{ij}=u_{ij}+ua_{ij}$. Let $T$ denote the unit tangent bundle of
$B_r(0)$ with respect to $g_0$. Assume the auxiliary function $z$
on $T$,  $z(x, e_p)=\rho^2\D^2 \wtt u (e_p, e_p)$,  attains its
maximum at $x_0$ and in direction $e_1 =(1, 0, \cdots, 0)$, where
$\rho(x)=(1-\frac{|x|^2}{r^2})^+$. In an orthonormal frame at
$x_0$, we may assume by a rotation of axes that $\{U_{ij}\}$ is
diagonal at $x_0$. Then at $x_0$, $F^{ij}$ is diagonal and
$$\align
 0 & =(\log z)_i  =\frac {2\rho_i}\rho
  +\frac {\wtt u_{11,i}} {\wtt u_{11}},\tag 3.16\\
 0& \ge (\log z)_{ii} =(\frac {2\rho_{ii}}{\rho}
              -\frac{6\rho_i^2}{\rho^2})
  +\frac {\wtt u_{11, ii}}{\wtt u_{11}}.\tag 3.17\\
 \endalign   $$
By (3.16), the gradient estimate, and the Ricci identities,
$$U_{ij, 11}=u_{ij 11}-\frac {u_{k1}^2}{u}\delta_{ij}
                +O(\frac{1+u_{11}}{\rho})
 = u_{11 ij}-\frac {u_{k1}^2}{u}\delta_{ij}
                +O(\frac{1+u_{11}}{\rho}).\tag 3.18$$
Hence we obtain
$$\align
 0 & \ge \sum_i F^{ii}(\log z)_{ii}
  \ge -\frac {C}{\rho^2}\Cal F
 + F^{ii}\frac {\wtt u_{11,ii}}{\wtt u_{11}}\\
   & \ge -\frac {C}{\rho^2}\Cal F
    +\frac {u_{11}^2}{2u \wtt u_{11}}\Cal F
    +\frac {1}{\wtt u_{11}} \D_k^2 \mu(f).\\
    \endalign  $$
Since $\mu(t)=t^{1/k}$, we have $\Cal F\ge C>0$. Hence (3.14)
holds. $\square$

The second order derivative estimate (3.14) was established in
[GW1]. As the proof is straightforward, we included it here for
completeness. The estimate is also similar to that in [GW4] for
the equation
$$\det \big(\D^2 u -\frac {|\D u|^2}{2u}I+\frac u2I\big)
       =f(x, u, \D u) \ \ \ \text{in}\ \ \Om\subset S^n\tag 3.19$$
which arises in the design of a reflector antenna, where $I$ is
the unit matrix. See also [W2] for $n=2$.

By Lemma 3.2, equation (3.1) becomes a uniformly elliptic
equation. By the Evans-Krylov estimates and linear theory [GT], we
have the following interior estimates.

\proclaim{Theorem 3.1} Let $u\in C^{3,1}$ be a positive solution
of (3.1) in a geodesic ball $B_r(0)\subset \M$. Suppose $f>0, \in
C^{1,1}$. Then for any $\alpha\in (0, 1)$,
$$\|u\|_{C^{3, \alpha} (B_{r/2}(0))}\le C , \tag 3.20$$
where $C$ depends only on $n, k, r$, $\inf_\M u$, and $g_0$.
\endproclaim

Theorem 3.1 also holds for equation (3.3) with $f=\kappa u^{-p}$
for a constant $p>0$ and a smooth, positive function $\kappa$.

\vskip25pt

{\bf 3.2. The parabolic equation}. It is more convenient to study
the parabolic equation for the function $w=\log u$. In this
section we will extend the a priori estimates in \S3.1 to the
parabolic equation
$$F[w]-w_t=\mu(f),\tag 3.21$$
where $F[w]=\mu[\sigma_k(\lam(W))]$, and
$$W=\D^2 w+\D w\otimes \D w -\frac 12 |\D w|^2g_0+A_{g_0}. $$
When $f=e^{-2kw}$, a stationary solution of (3.21) satisfies the
equation
$$\sigma_k(\lam(W))=e^{-2kw},$$
which is equivalent to (3.1).

We choose a monotone increasing function $\mu$ such that $F$ is
concave in $D^2 w$ and
$$\mu(t)=\cases
 t^{1/k}\ \ \ & t\ge 10,\\
 \log t\ \ \ &t\in (0, \frac 1{10}),\\
 \endcases $$
and furthermore
$$(t-s)(\mu(t)-\mu(s))\ge c_0(t-s)(t^{1/k}-s^{1/k}) \tag 3.22$$
for some constant $c_0>0$ independent of $t$. Condition (3.22)
will be used in the next section.

We say $w$ is $k$-admissible if for any fixed $t$, $w$ is
$k$-admissible as a function of $x$. Denote $Q_r=B_r(0)\times (0,
r^2]$. In the following lemmas we establish interior (in both time
and spatial variables) a priori estimates for $w$,

\proclaim {Lemma 3.3} Let $w$ be a $k$-admissible solution of
(3.21) on $Q_r$. Then we have the estimates
$$|\D_x w(0, r^2)| \le C, \tag 3.23 $$
where $C$ is independent of $\sup w$, if $f=\kappa (x) e^{-pw}$
for some constant $p>0$ and smooth, positive function $\kappa$.
\endproclaim

\noo{\it Proof}. The proof is similar to that of Lemma 3.1. We
outline the proof here. Let $u=e^w$. Then $u$ satisfies the
equation
$$\wtt F[u]-\frac {u_t}{u}=\mu(f),\tag 3.24$$
where
$$\wtt F[u]=\mu\big[\frac {1}{u^k}\sigma_k
      \big(\lam(\D^2 u-\frac {|\D u|^2}{2u}g_0+uA_{g_0})\big)\big].$$
Let $z=\big(\frac {|\D u|}{u}\big)^2\rho^2$ be the auxiliary
function as in the proof of Lemma 3.1. Here we choose
$$\rho(x, t)=\frac{t}{r^2}(1-\frac{|x|^2}{r^2})^+. $$
Suppose $z$ attains its maximum at $(x_0, t_0)$. Then $t_0>0$. By
a rotation of axes we assume $|\D u|=u_1$. Then at $(x_0, t_0)$,
$z_i=0$, $\{z_{ij}\}\le 0$, and $z_t\ge 0$. Hence we have (3.4),
(3.5) and
$$\frac {u_{1t}}{u}-\frac {u_1u_t}{u^2}
                 +\frac{\rho_t}{\rho}\ge 0 .\tag 3.25$$
Differentiating equation (3.24) we obtain (3.6) with $F^{ij}$ and
$\Delta$ replaced by
$$\align
\wtt F^{ij}(r)&
 =\frac {\p}{\p r_{ij}} \mu[\frac {1}{u^k}\sigma_k(\lam(r))]
 =\frac{\mu'}{u^k}\frac{\p}{\p r_{ij}}\sigma_k(\lam(r)),\\
\Delta & =[\frac {u_{1t}}{u}-\frac {u_1u_t}{u^2}]
           +\frac {ku_1\mu'}{u^{k+1}}\sigma_k(\lam)
           +[\D_1 \mu(f)
           -\wtt F^{ij}\D_1 (a_{ij}u)]\\
 &\ge -\frac{\rho_t}{\rho}+\D_1 \mu(f)
           -\wtt F^{ij}\D_1 (a_{ij}u).\\
 \endalign $$
By Lemma 2.3 and our choice of $\mu$, $\wtt\F =\sum \wtt F^{ii}$
has a positive lower bound,
$$\wtt\F\ge \frac {\mu'(u^{-k}\sigma_k(\lam))}{u^k}
            \sigma_k^{(k-1)/k} (\lam(U))
            \ge \frac Cu$$
for some $C>0$ depending only on $n, k$. Hence similarly as the
proof of Lemma 3.1 we have (3.7). From (3.8), we obtain estimate
(3.23). $\square$

Since the constant $C$ in (3.23) is independent of $\sup w$, hence
we have a similar Harnack type inequality as Corollary 3.2.

\proclaim {Lemma 3.4} Let $w$ be a $k$-admissible solution of
(3.21) on $Q_r$.  Then we have the estimate
$$|\D_x^2 w(0,r^2)| \le C, \tag 3.26 $$
where $C$ depends only on $n, k$, $r, \mu$, $\inf w$, $\sup w$,
and $\|A_{g_0}\|_{C^2}$.
\endproclaim

\noo{\it Proof}. Differentiating equation (3.21) twice, we get
$$\align
 F^{ij} W_{ij, k} & =w_{tk}+\D_k \mu(f), \tag 3.27\\
 F^{ij} W_{ij, kk} &=-F^{ij, rs}W_{ij, k}W_{rs, k}+w_{tkk}
                 +\D_k^2 \mu(f)
       \ge w_{tkk}+\D_k^2 \mu(f).\tag 3.28\\
 \endalign $$
where $F^{ij}=\frac{\p F}{\p W_{ij}}$, $W_{ij, k}=\D_k W_{ij}$,
and $F^{ij, rs}=\frac {\p^2 \mu(\sigma_k(\lam(W)))}{\p W_{ij}\p
W_{rs}}$. Denote $\wtt w_{ij}=w_{ij}+a_{ij}$,
$a_{ij}=(A_{g_0})_{ij}$. Let $T$ denote the unit tangent bundle of
$\M$ with respect to $g_0$. Consider the auxiliary function $z$
defined on $T\times [0, r^2]$, given by $z=\rho^2\big(\D^2\wtt
w+(\D w)^2\big)(e_p, e_p)$, where $\rho$ is the cut-off function
in the proof of Lemma 3.3. Assume that $z$ attains its maximum at
$(x_0, t_0)$ and in direction $e_1=(1, 0, \cdots, 0)$. We choose
an orthonormal frame at $(x_0, t_0)$, such that after a rotation
of axes, $\{W_{ij}\}$ is diagonal. Then $F^{ij}$ is diagonal and
at $(x_0, t_0)$,
$$\align
 0 & =(\log z)_i  =\frac {2\rho_i}\rho
  +\frac {\wtt w_{11,i}+2w_1w_{1i}} {\wtt w_{11}+w_1^2},\tag 3.29\\
 0 & \le (\log z)_t  =\frac {2\rho_t}\rho
  +\frac {w_{11t}+2w_1w_{1t}} {\wtt w_{11}+w_1^2},\tag 3.30\\
 0& \ge (\log z)_{ii} =(\frac {2\rho_{ii}}{\rho}
              -\frac{6\rho_i^2}{\rho^2})
  +\frac {\wtt w_{11, ii}+2w_1w_{1ii}+2w_{1i}^2}
                   {\wtt w_{11}+w_1^2}.\tag 3.31\\
 \endalign   $$
We have, by (3.29) and the Ricci identities,
$$\align
W_{ij, 11} & =w_{ij 11}+w_{i11}w_j+w_{j11}w_i+2w_{i1}w_{j1}
    -w_{k1}^2\delta_{ij}+O(\frac 1\rho (\wtt w_{11}+w_1^2))\\
   & = w_{11 ij}+2w_{i1}w_{j1} -w_{k1}^2\delta_{ij}
              +O(\frac 1\rho (\wtt w_{11}+w_1^2)). \\
     \endalign $$
Hence we obtain
$$\align
 0 & \ge \sum_i F^{ii}(\log z)_{ii} -(\log z)_t\\
  &\ge -\frac {C}{\rho^2}\Cal F +\frac{1}{\wtt w_{11}+w_1^2}
   F^{ii} (\wtt w_{11, ii}+2w_1w_{ii1}+2w_{i1}^2)
   -\frac {w_{11t}+2w_1w_{1t}}{\wtt w_{11}+w_1^2}
     -\frac {2\rho_t}\rho\\
  &\ge -\frac {C}{\rho^2}\Cal F +\frac{1}{\wtt w_{11}+w_1^2}
   F^{ii}[ (W_{ii, 11}+w_{k1}^2)+2w_1w_{1ii}]
   -\frac {w_{11t}+2w_1w_{1t}}{\wtt w_{11}+w_1^2}
       -\frac {2\rho_t}\rho\\
&\ge -\frac {C}{\rho^2}\Cal F +\frac{1}{\wtt w_{11}+w_1^2}
   ( F^{ii}W_{ii, 11}- w_{11t})+w_{11}\Cal F
   +\frac {2w_1}{\wtt w_{11}+w_1^2}(F^{ii}w_{ii1}-w_{t1})
          -\frac {2\rho_t}\rho\\
&\ge -\frac {C}{\rho^2}\Cal F
 +\frac{1}{\wtt w_{11}+w_1^2}\D_1^2 \mu(f)+
 w_{11}\Cal F+\frac {2w_1}{\wtt w_{11}+w_1^2}
   \D_1 \mu(f) -\frac {2\rho_t}\rho.\tag 3.32\\
   \endalign $$
By our choice of $\mu$, $\F\ge C$ for some $C$ depending only on
$n, k$. We obtain $w_{11}\rho^2\le C$ at $(x_0, t_0)$. Whence
$z(0, r^2)\le z(x_0, t_0)\le C$.  $\square$

\noo{\bf Remark 3.1}.  The a priori estimates (3.23) and (3.26)
also hold for the equation
$$\frac 1a \mu(a^k\sigma(\lam(W))-w_t=\frac 1a \mu(a^kf),
                                           \tag 3.33 $$
where $a>0$ is a constant. We claim that the constants $C$ in
(3.23) and (3.26) are independent of $a\ge 1$.

First we note that this is obvious for (3.23) as the gradient
estimate is independent of the choice of $\mu$.

For the estimate (3.26), we have, by Lemma 2.3,
$$\align
\F & =\sum_i \frac {\p}{\p W_{ii}}
           \bigg[\frac 1a\mu (a^k\sigma_k(\lam(W))\bigg]\\
   & =(n-k+1)a^{k-1}\sigma_{k-1}(\lam(W))\mu'\\
   &\ge Ca^{k-1}\sigma_k^{(k-1)/k} \mu'(a^k\sigma_k)\\
   &\ge C\inf_{t>0}t^{(k-1)/k}\mu'(t),\\
   \endalign $$
By our choice of $\mu$, $\inf_{t>0}t^{(k-1)/k}\mu'(t)\ge C>0$.
Hence $\F>C>0$.

Therefore by (3.32) it suffices to show that
$$|\D_1 g|+|\D_1^2 g|\le C$$
for some $C>0$ independent of $a\ge 1$, where $g=\frac 1a
\mu(a^kf)$. By our choice of $\mu$,
$$g=\cases
 \mu(f)\ \ \ &\text{if}\ a^kf>10\\
 \frac 1a(k\log a+\log f)\ \ &\text{if}\ a^kf<\frac 1{10},\\
 \endcases$$
Hence $\sup (|\D_1g|+|\D_1^2 g|)$ is independent of $a\ge 1$  if
$a^kf>10$ or $a^kf<\frac 1{10}$. When $a^kf \in (\frac 1{10},
10)$, we can choose $\mu$ properly such that $\sup |\D_1^2 g|$ is
independent of $a\ge 1$. Alternatively one can compute directly
$$\align
|\D_1g| & =|a^{k-1}f^{(k-1)/k}\mu'(a^k f)|\,|k\D_1f^{1/k}|\\
       &\le C|\D_1f^{1/k}| ,\\
|\D_1^2g| & =|a^{k-1}\mu'(a^k f)\D_1^2 f|
          +\big|a^{k-1}\mu'(a^k f)\frac {(\D_1f)^2}{f}\big|\,
          \big|\frac {a^k f\mu''(a^k f)}{\mu'(a^k f)}\big|\\
  &\le C\frac{|\D^2_1 f|}{f^{1-\frac 1k}}|
       + C\frac {|\D_1 f|^2}{f^{2-\frac 1k}},\\
\endalign $$
where $C$ depends on $\sup_{t>0}t^{1-1/k}\mu'(t)$ and
$\sup_{t>0}\frac {t\mu''(t)}{\mu'(t)}<\infty$. Hence the estimate
(3.26) is independent of $a\ge 1$.

Remark 3.1 will be used in the next section. The choice of
function $\mu$ in the parabolic equation (3.21) is critical.

\proclaim {Lemma 3.5} Let $w$ be a $k$-admissible solution of
(3.21) on $Q_r$.  Then we have the estimates
$$| w_t(0,r^2)| \le C, \tag 3.34$$
where $C$ depends only on $n, k$, $r, \mu$, $\inf w$, $\sup w$,
and $\|A_{g_0}\|_{C^2}$.
\endproclaim

\noo{\it Proof}. From the equation (3.21) and by the estimate
(3.26) we have an upper bound for $w_t$. It suffices to show that
$w_t$ is bounded from below. Let
$z=\frac{w_t}{(M-w)^\alpha}\rho^\beta$, where $M=2\sup_{Q_r} |w|$,
$\rho$ is the cut-off function as above. Suppose $\min_{Q_r} z$
attains its minimum at $(x_0, t_0)$, $t_0>0$.  Then at the point
we have $z_t\le 0$, $z_i=0$ and the matrix $\{z_{ij}\}\ge 0$,
namely
$$\align
& \frac{w_{tt}}{w_t}+\alpha\frac {w_t}{M-w}
                      +\beta \frac{\rho_t}{\rho}\ge 0,\tag 3.35\\
& \frac{w_{ti}}{w_t}+\alpha\frac {w_i}{M-w}
          +\beta \frac{\rho_i}{\rho}=0\ \ \ i=1, \cdots, n,\tag 3.36\\
& \{ \frac{w_{ijt}}{w_t} -\frac{w_{it}w_{jt}}{w_t^2}
      +\alpha\frac {w_{ij}} {M-w} +\alpha\frac {w_iw_j}{(M-w)^2}
         +\beta\frac{\rho_{ij}}{\rho}
            -\beta\frac {\rho_i\rho_i}{\rho^2}  \}\le 0,\tag 3.37\\
\endalign$$
where we have changed the direction of the inequalities as we
assume that $w_t<0$. Differentiating equation (3.21) gives
$$F^{ij}W_{ijt}-w_{tt}=\frac{\p}{\p t} \mu(f). \tag 3.38$$
Hence by (3.35),
$$\align
\alpha \frac {w_t}{M-w}
    &\ge -\frac{w_{tt}}{w_t}-\beta\frac{\rho_t}{\rho}\\
        & =\frac{-1}{w_t}F^{ij}W_{ijt}+
        \frac{1}{w_t}\frac {\p}{\p t }\mu(f)
          -\beta\frac{\rho_t}{\rho}\\
        \endalign $$
By (3.36), the matrix in (3.37) is equal to
$$\{\frac{w_{ijt}}{w_t}  +\frac {\alpha w_{ij}} {M-w}
          +\frac {\alpha(1-\alpha)w_iw_j}{(M-w)^2}
            -\frac{2\alpha\beta w_i\rho_j}{(M-w)\rho}
              +\beta\frac{\rho_{ij}}{\rho}
               -\beta(1+\beta)\frac {\rho_i\rho_i}{\rho^2} \}\le 0. $$
We have
$$\align
\frac{-1}{w_t}F^{ij}W_{ijt}
 &=\frac{-1}{w_t} F^{ij}
          (w_{ijt}+w_{it}w_j+w_{jt}w_i-w_kw_{kt}\delta_{ij})\\
 &\ge F^{ij}(\frac {\alpha w_{ij}} {M-w}
          +\frac {\alpha(1-\alpha)w_iw_j}{(M-w)^2}
            -\frac{2\alpha\beta w_i\rho_j}{(M-w)\rho}
              +\beta\frac{\rho_{ij}}{\rho}
               -\beta(1+\beta)\frac {\rho_i\rho_i}{\rho^2})\\
  &\ \ \ +F^{ij}(2\alpha\frac{w_iw_j}{M-w}
                 +2\beta\frac{\rho_iw_j}{\rho}
   -\alpha\frac{|\D w|^2}{M-w}\delta_{ij}
               -\beta\frac{w_k\rho_k}{\rho}\delta_{ij}\}\\
  &\ge \frac{\alpha}{M-w}F^{ij}(w_{ij}+2w_iw_j
                         -|\D w|^2\delta_{ij})\\
   &\ \ \ +\frac{\alpha}{M-w}F^{ij}(\frac{(1-\alpha)w_iw_j}{M-w}
             -2\beta\frac{w_i\rho_j}{\rho})-\frac{C}{\rho^2}, \\
                 \endalign $$
where the constant $C$ depends on the gradient estimate (3.23) and
the second derivative estimate (3.26). Choose $\alpha=\frac 12$.
By the Holder inequality,
$$F^{ij}(\frac{w_iw_j}{2(M-w)}-2\beta\frac{w_i\rho_j}{\rho})\ge-\frac{C}{\rho^2}.$$
By the $k$-admissibility, $F^{ij}W_{ij}\ge 0$. Hence we obtain
$$\frac{-1}{w_t} F^{ij}W_{ijt}
 \ge \frac {\alpha} {M-w} F^{ij} (W_{ij}+w_iw_j
           -\frac 12|\D w|^2\delta_{ij}-a_{ij}) -\frac{C}{\rho^2}
 \ge -\frac{C}{\rho^2}. $$
 It follows that
$$\alpha \frac {w_t}{M-w}\ge -\frac{C}{\rho^2}
             +\frac {1}{w_t}\frac {\p}{\p t }\mu(f)
                  -\beta\frac{\rho_t}{\rho}. $$
Now we choose $\beta=2$. Then we obtain
$$z(x_0, t_0)=\frac{w_t}{(M-w)^{1/2}}\rho^2(x_0, t_0)\ge -C. $$
It follows that $z(0, r^2)\ge z(x_0, t_0)\ge -C$. Hence $w_t$ is
bounded from below. $\square$

\proclaim{Theorem 3.2} For any $w_0\in\Phi_k$, there is a smooth
$k$-admissible solution $w\in C^{3, 2}(\M\times [0, T))$ of (3.21)
with $w(\cdot, 0)=w_0$ on a maximal time interval $[0, T)$. If
$T<\infty$, we have $\inf_\M w(\cdot, t)\to -\infty$ as $t\nearrow
T$.
\endproclaim

\noo{\it Proof}. First we point out that a $k$-admissible solution
of (3.21) is locally bounded. Indeed, at the minimum point of $w$,
by equation (3.21) we have
$$w_t=F[w]-\mu(f)
     \ge \mu\big(\sigma_k(\lam(A_{g_0}))\big)-\mu(f).$$
Hence locally in time the solution is bounded from below. By the
interior gradient estimate (3.23), the solution is also bounded
from above. Therefore by Lemmas 3.3-3.5, equation (3.21) is
uniformly parabolic. By Krylov's regularity theory, we obtain the
$C^{3, 2}$ a priori estimate for (3.21), and so the local
existence follows. Let $[0, T)$ be the maximal time interval for
the solution. If $T<\infty$, we must have $\inf_\M w(\cdot, t)\to
-\infty$ as $t\nearrow T$. $\square$

\vskip10pt

\noo {\bf Remark 3.2}. The a priori estimates in \S 3.1 and \S 3.2
can be extended to the quotient equation
$$\sigma_{k, l}\big(\lam(\D^2 u
    -\frac {|\D u|^2}{2u}g_0+uA)\big)=u^{l-k}\ \ \ \
    (1\le l<k\le n) \tag 3.39$$
and its parabolic counterpart,  where $\sigma_{k, l}(\lam)=\frac
{\sigma_k}{\sigma_l}(\lam)$. Indeed, let $\mu$ be a monotone
increasing function such that $\mu[\sigma_{k, l}(\lam)]$ is
concave. Write equation (3.39) in the form $F[u]=\mu(f)$ as (3.3).
Then the proof for the second derivative estimates (3.14) and
(3.26) can be extended to the quotient equation without change.
For the gradient estimates (3.2) and (3.23), denote
 $\sigma_{k-1; i}(\lam)=\frac {\p}{\p\lam_i}\sigma_k(\lam)$.
By Newton's inequality [LT],
$$\align
\frac {\p\sigma_{k, l}}{\p \lam_i}(\lam)
 & =\frac{\sigma_l\sigma_{k-1;i}
                -\sigma_k\sigma_{l-1;i}}{\sigma_l^2}\\
 & =\frac {\sigma_{l;i}\sigma_{k-1;i}
                -\sigma_{k;i}\sigma_{l-1;i}}{\sigma_l^2}
 \ge \frac {n(k-l)}{k(n-l)}
          \frac{\sigma_{l;i}\sigma_{k-1;i}}{\sigma_l^2}.\tag 3.40\\
          \endalign$$
As before we arrange the eigenvalues in the descending order
$\lam_1\ge\cdots\ge\lam_n$. Then by Lemma 2.3(v),
$\sigma_{l;i}(\lam)\ge C\sigma_l(\lam)$ when $i\ge l+1$. Hence
$\frac {\p\sigma_{k, l}}{\p \lam_k}(\lam)
 \ge C \sum_i \frac {\p\sigma_{k, l}}{\p \lam_i}(\lam)$,
namely $F^{kk}\ge C\Cal F$. Hence the proof of (3.2) and (3.23)
can also be extended to the quotient equation. But we need to
replace (3.11) in Case 2 by
$$\align
A & \ge \frac 12\mu'
    \frac {\p \sigma_{k, l}}{\p \lam_{k-1}}\wtt u_{k-1\,k-1}^2
  \ge \frac 12\mu' \frac {n(k-l)}{k(n-l)} \frac
   {\sigma_{l;k-1}\sigma_{k-1;k-1}}{\sigma_l^2}\wtt u_{k-1\,k-1}^2\\
 &  \ge C\frac{\sigma_{l;k-1}\sigma_{k-1;k-1}}
                                 {\sigma_l^2}\lam_{k-1}^2
   \ge C\frac{\sigma_{k-1;k-1}}{\sigma_l}\lam_{k-1}^2
  \ge Cb^2\Cal F ,\\
  \endalign $$
and (3.12) in Case 3 by
$$A\ge \frac 12\Sigma_{i\ne i^*} F^{ii}\wtt u^2_{ii}
   \ge \frac 12 F^{ll}\wtt u^2_{ll}
   \ge C\frac {\sigma_{n-1, l}}{\sigma_l}\wtt u_{ll}^2
   \ge C\lam_l\lam_n \frac {\sigma_{n-1, n}}{\sigma_l}
   \ge Cb^2\Cal F,$$
where we have used $F^{ll}\ge C\sigma_{k-1, l}/\sigma_l$ by
(3.40).

For the corresponding parabolic equation (3.21), where
$F[w]=\mu\big(\sigma_{k, l}(\lam(W))\big)$,  choose a monotone
increasing function $\mu$ such that $F$ is concave in $\lam$, and
$\mu(t)=t^{1/(k-l)}$ when $t\ge 10$,  $\mu(t)=\log t$ for $t>0$
small. Then we can prove Theorem 3.2 for the quotient equation in
the same way as before.

We remark that the a priori estimates for (3.39), and for its
parabolic counterpart on locally conformally flat manifolds, were
obtained in [GW3].

\vskip5pt

We also note that the a priori estimates in \S 3.1 can be extended
to the more general equation
$$s_k=\sum_{l=0}^{k-1}\beta_ls_l,\tag 3.41$$
where $\beta_l$ are nonnegative constants, $\sum_l \beta_l>0$, and
$s_k=u^k\sigma_k(\lam(\D^2 u-\frac {|\D u|^2}{2u}+uA))$ is the
$\kk$-curvature.

\vskip25pt

{\bf 3.3. Counterexamples}. Theorem 3.1 applies to solutions of
(3.1) with eigenvalues in the positive cone $\Ga_k$. The a priori
estimate (3.14) relies critically on the negative sign of the term
$\frac {|\D u|^2}{2u}$, which yields the dominating term
$u_{k1}^2$ in (3.18). Equation (3.1) has another elliptic branch,
namely when the eigenvalues $\lam$ lie in the negative cone
$-\Ga_k$. An open problem is whether the a priori estimate (3.14)
holds for solutions with eigenvalues in the negative cone
$-\Ga_k$. This is also an open problem for equations from optimal
transportation [MTW], in particular the reflector antenna design
problem (3.19). Here we give a counter example to the regularity.
Our example is a modification of the Heinz-Levy counterexample in
[Sc].

We will consider the two dimensional case. By making the change
$u\to -u$, we consider equation
$$\det (u_{ij}+|\D u|^2I+a_{ij})=f\tag 3.42$$
with positive sign before the term $|\D u|^2$, where $f$ is a
$C^{1,1}$ positive function to be determined. We want to show that
equation (3.42) has no interior a priori estimates for solutions
with eigenvalues in the positive cone.

Set
$$u(x)=\frac {b}{2}x_2^2+\phi(x_1),\tag 3.43$$
where $b$ is constant, $\phi$ is an even function. Let
$$a_{11} = -b^2x_2^2,\ \ \ \ a_{12} = 0,
       \ \ \ \ a_{22} = -b-b^2x_2^2. \tag 3.44$$
Then equation (3.42) becomes
$$(\phi^{''}+{\phi'}^2){\phi'}^2=f. \tag 3.45$$
Let $\psi=(\phi')^3$. Then $\psi$ satisfies the equation
$$\frac 13\psi'+\psi^{4/3}=f.\tag 3.46 $$
Let
$$\psi(x_1)  =x_1-\frac 97x_1^{7/3} .\tag 3.47$$
Then
$$f(x)=\frac 13-x_1^{4/3}+(x_1-\frac 97 x_1^{7/3})^{4/3}\tag 3.48$$
is a positive $C^2$ function, but the solution $u\not\in C^2$.

If instead of (3.44), we choose
$$a_{11} = c_0-b^2x_2^2,\ \ \ \ a_{12} = 0,
       \ \ \ \ a_{22} = \eps-b-b^2x_2^2, \tag 3.49$$
where $c_0,\eps $ are constants, $\eps>0$ small. Then we have the
equation
$$(\phi^{''}+{\phi'}^2+c_0)(\eps+{\phi'}^2)=f. \tag 3.50$$
Let $f\equiv 1$ and denote $g=\phi'$. Then $g(0)=0$ and $g$
satisfies
$$g'=\frac 1{\eps+g^2}-g^2-c_0.\tag 3.51 $$
This equation has a unique solution $g_\eps$. Obviously the gradient
of $g_\eps$ is not uniformly bounded. Hence there is no interior
$C^{1,1}$ a priori estimate for equation (3.43). Note that the
matrix $A=(a_{ij})$ can either be in the positive cone or in
negative cone by choosing proper constants $b, c_0$.

Write equation (3.1) in the form
$$\sigma_k(\lam(\D^2 w-\D w\otimes \D w+\frac 12 |\D w|^2I+A))
                      =f.\tag 3.52$$
Then similarly as above we can construct a sequence of functions
satisfying equation (3.52) with $f=1$ whose second derivatives are
not uniformly bounded.

\vskip10pt

\noo{\bf Remark 3.3}. In many situations [MTW, W2] there arise
equations of the form
$$\sigma_k(\lam(D^2 u+A(x, u, Du))=f,\tag 3.53$$
where $A$ is a matrix. From the discussions in this section, we
see that the interior a priori estimates hold in general when $A$
is negative definite with respect to $Du$, and do not hold if $A$
is positive definite. When $A=0$, there is no interior regularity
in general, but if the solution vanishes on the boundary, interior
a priori estimates have been established in [CW2].

\newpage

\vskip30pt

\centerline{\bf 4. Proof of Lemma 2.1}

\vskip10pt

{\bf 4.1. Existence of solutions in the sub-critical growth case}.
In this subsection we first study the existence of $k$-admissible
solutions, for $2\le k<\frac n2$, to equation (2.5) in the
subcritical growth case $1<p<\frac{n+2}{n-2}$. We then extend the
existence result to the critical case $p=\frac {n+2}{n-2}$ in \S
4.2 by the blow-up argument. In \S 4.3 we consider the case
$k=\frac n2$.

\proclaim{Theorem 4.1} Suppose $2\le k<\frac n2$. Then for any
given $1<p<\frac {n+2}{n-2}$, there is a solution $v_p$ of (2.5)
with $J_p(v_p)=c_p>0$, where $J_p, c_p$ are defined respectively
in (2.6) and (2.8). Moreover, the set of solutions of (2.5) is
compact.
\endproclaim

A solution of (2.5) is a critical point of the functional $J=J_p$.
To study the critical points of the functional $J$, we will employ
the parabolic equation (3.21), which can also be written in terms
of $v$ as (ignoring a coefficient $\frac {2}{n-2}$ before $v_t$)
$$F[v]+\frac{v_t}{v}=\mu(f(v)), \tag 4.1$$
where $f(v)=v^{\frac{4k}{n-2}-\eps}$,
$F[v]=\mu(\sigma_k(\lam(\frac Vv )))$, $\mu$ is the function in
(3.21), and
$$\eps=\frac {n+2}{n-2}-p.\tag 4.2$$

Write functional (2.6) in the form
$$J(v)=\frac {n-2}{2n-4k}\int_{(\M, g_0)}
       v^{\frac {2n-4k}{n-2}}\sigma_k(\lam(\frac Vv))
   -\frac {1}{p+1}\int_{(\M, g_0)}
         v^{\frac {2n-4k}{n-2}}v^{\frac {4k}{n-2}-\eps}. \tag 4.3$$
Equation (4.1) is a descent gradient flow of the functional $J$,
$$\align
\frac {d}{dt}J(v)
  &= \int_{(\M, g_0)}v^{\frac {2n-4k}{n-2}-1}
  \big[\sigma_k(\lam(\frac Vv))-v^{\frac{4k}{n-2}-\eps}\big]v_t\\
  &= -\int_{(\M, g_0)}v^{\frac {2n-4k}{n-2}}
   \big[\sigma_k(\lam(\frac Vv))-v^{\frac{4k}{n-2}-\eps}\big]
   \, \big[\mu\big(\sigma_k(\lam(\frac Vv))\big)
                   -\mu\big(v^{\frac{4k}{n-2}-\eps}\big)\big]
        \le 0.\tag 4.4 \\
        \endalign $$

\vskip5pt

Given an initial $k$-admissible function $v_0$, by Theorem 3.2,
the flow (4.1) has a unique smooth positive solution $v$ on a
maximal time interval $[0, T)$, where $T\le \infty$.

\proclaim{Lemma 4.1} Suppose $J(v(\cdot, t))$ is bounded from
below for all $t\in (0, T)$. If $v(\cdot, t)$ is uniformly
bounded, then either $v(\cdot, t)\to 0$ or there is a sequence
$t_j\to\infty$ such that $v(\cdot, t_j)$ converges to a solution
of (2.5).
\endproclaim

\noo{\it Proof}. By the assumption that $v(\cdot, t)$ is uniformly
bounded, we have $T=\infty$.  At the maximum point of $v(\cdot,
t)$, by equation (4.1) we have
$$v_t\le v[\mu(f(v))-\mu(\sigma_k(\lam(A_{g_0})))] . \tag 4.5$$
Hence if $\sup v(\cdot, t_0)$ is sufficiently small at some $t_0$,
by the assumptions $ g_0 \in\Ga_k$ and $v>0$, we have $v(\cdot,
t)\to 0$ uniformly. Therefore if $v$ does not converges to zero
uniformly, by the gradient estimate (3.23), we have $v\ge c$ for
some constant $c>0$. In the latter case, by Theorem 3.2 and the
assumption that $v$ is uniformly bounded,  we have $\|v(\cdot,
t)\|_{C^3(\M)}\le C$ for any $t\ge 0$.

Choose a sequence $t_j\to\infty$ such that
$$\frac {d}{dt}J(v(\cdot, t_j))\to 0.\tag 4.6$$
By the above $C^3$ a priori estimate, we may abstract a
subsequence, still denoted as $t_j$, such that $v(\cdot, t_j)$
converges in $C^{2, \alpha}$. By (4.4) we concludes that $v(x,
t_j)$ converges as $j\to\infty$ to a solution of (2.5). $\square$

\proclaim{Lemma 4.2} Suppose $J(v(\cdot, t))$ is bounded from
below for all $t\in (0, T)$. Then $T=\infty$ and $v(\cdot, t)$ is
uniformly bounded.
\endproclaim

\noo{\it Proof}. Suppose to the contrary that there exists a
sequence $t_j\nearrow T$ such that $m_j=\sup v(\cdot,
t_j)\to\infty$. Assume the maximum is attained at $z_j\in\M$. By
choosing a normal coordinate centered at $z_j$, we may identify a
neighbourhood of $z_j$ in $\M$ with the unit ball in $\R^n$ such
that $z_j$ becomes the origin. We make the local transformation
$$\align
v_j (y, s) & =m_j^{-1}v(x, t),\tag 4.7\\
y & = m_j^{\frac{2}{n-2}-\frac{\eps}{2k} } x, \\
s & = m_j^{\frac {4}{n-2}-\frac\eps k} (t-t_j).\\
\endalign $$
For the transformation $x\to y$, more precisely it should be
understood as a dilation of $\M$, regarded as a submanifold in
$\R^N$ for some $N>n$ with induced metric. Denote $\M_j=\{Y=
m_j^{\frac{2}{n-2}-\frac{\eps}{2k}} X\ |\ X\in\M\subset\R^N\}$,
with induced metric from $\R^N$. Then we have $0<v_j(y, 0)\le
m_j^{-1}v(0,t_j)=1$, $v_j$ is defined for $y\in\M_j$ and $s\le
s_0$, where by (4.5), $s_0>0$ is a positive constant independent
of $j$. Moreover $v_j$ satisfies the equation
$$m_j^{-\frac{4}{n-2}+\frac{\eps}{k}}
  \mu\big[m_j^{\frac {4k}{n-2}-\eps}
               \sigma_k(\lam(\frac {V_j}{ v_j}))\big]
                           + \frac {(v_j)_s}{v_j}
   =m_j^{-\frac{4}{n-2}+\frac{\eps}{k}}
        \mu (m_j^{\frac {4k}{n-2}-\eps}f(v_j)). \tag 4.8 $$

By direct computation,
$$\align
\int_{\M_j} {v_j}^{\frac {2n-4k}{n-2}}
                         \sigma_k(\lam(\frac{ V_j}{v_j}))dy
 &= m_j^{\eps(1-\frac{n}{2k})} \int_{\M} {v}^{\frac {2n-4k}{n-2}}
                         \sigma_k(\lam(\frac{V}{v}))dx \tag 4.9\\
\int_{\M_j}{v_j}^{\frac {2n}{n-2}-\eps} dy
 &=m_j^{\eps(1-\frac{n}{2k})}\int_{\M}{v}^{\frac {2n}{n-2}-\eps}dx\\
 \endalign $$
Hence
$$\align
J(v_j, \M_j)
 & =:\frac {n-2}{2n-4k}\int_{\M_j} {v_j}^{\frac {2n-4k}{n-2}}
                         \sigma_k(\lam(\frac{ V_j}{v_j}))dy
   -\frac {1}{p+1}\int_{\M_j}{v_j}^{\frac {2n}{n-2}-\eps}dy\\
 & =m_j^{\eps(1-\frac{n}{2k})} J(v, \M)\le C. \tag 4.10\\
 \endalign $$
We may choose $s_j\in (0, \frac 12s_0)$ such that
$$\frac{d}{ds} J(v_j (\cdot, s_j))\to 0. \tag 4.11$$
By (4.4), (4.11) is equivalent to
$$\align
\int_{\M_j} & {v_j}^{\frac {2n-4k}{n-2}}
    \bigg\{\sigma_k(\lam(\frac {V_j}{v_j}))
               -{v_j}^{\frac{4k}{n-2}-\eps}\bigg\}\cdot\\
 & \bigg\{m_j^{-\frac{4}{n-2}+\frac\eps k} \bigg[
 \mu\big(m_j^{\frac{4k}{n-2}-\eps}\sigma_k(\lam(\frac {V_j}{v_j}))\big)
 -\mu\big(m_j^{\frac{4k}{n-2}-\eps}v_j^{\frac{4k}{n-2}-\eps}\big)\bigg]
  \bigg\}\to 0.\\
  \endalign $$
By (3.22), we obtain
$$\int_{\M_j} {v_j}^{\frac {2n-4k}{n-2}}
   \bigg\{\sigma_k(\lam(\frac {V_j}{v_j}))
               -{v_j}^{\frac{4k}{n-2}-\eps}\bigg\}\cdot
\bigg\{ \big(\sigma_k(\lam(\frac {V_j}{v_j}))\big)^{1/k}
    -\big(v_j^{\frac{4k}{n-2}-\eps}\big)^{1/k}\bigg\}\to 0.\tag 4.12 $$
By the gradient estimate (3.23), $v_j+\frac {1}{v_j}$ (at $s=s_j$)
is locally uniformly bounded. Hence
$$\sigma_k(\lam(\frac {V_j}{v_j}))
          -v_j^{\frac{4k}{n-2}-\eps}\to 0\ \
           \text{in}\ \ L^{(k+1)/k}.\tag 4.13$$
Note that by Remark 3.1, $v_j$ are locally uniformly bounded in
$C^{1,1}$ and the convergence in (4.13) is locally uniform.

By extracting a subsequence we can assume that $v_j(\cdot, s_j)$
converges to a function $v_0\in C^{1,1}(\R^n)$ with $v_0(0)=1$. We
claim that $v_0$ is a smooth solution of the equation
$$F_0[v] := \sigma_k^{1/k}(\lam(V)) =v^p\tag 4.14$$
in $\R^n$. On the other hand, by the Liouville Theorem in [LL2],
there is no entire positive solution to (4.14) when $\eps>0$. This
is a contradiction. Hence Lemma 4.2 holds.

To prove that $v_0$ is a smooth solution of (4.14), we note that
since $v_0\in C^{1,1}(\R^n)$, $v_0$ is twice differentiable almost
everywhere. Suppose now that $F_0[v_0]>v_0^p$ at some point $x_0$
where $v_0$ is twice differentiable. Without loss of generality we
assume that $x=0$. Let
$$\phi(x)= v_0(0)+Dv_0(0)x+\frac 12 D_{ij}v_0(0)x_ix_j
                          +\frac \eps 2 |x|^2-\delta, $$
where $\eps, \delta$ are positive constants. By choosing $\delta$
sufficiently small we have
$$\phi>v_0\ \ \ \text{on}\ \ \p B_r(0)
 \ \ \ \ \text{and}\ \ \ \phi(0)<v_0(0). $$
Since $v_j\to v_0$ uniformly, we have $\phi>v_j$ on $\p B_r(0)$
and $\phi(0)<v_j(0)$ when $j$ is sufficiently large. Since $v_0$
is locally uniformly bounded in $C^{1,1}$, by the inequality
$F_0[v_0]>v_0^p$ we have $\lam(V_{v_0}-\eps I)\in \Ga_k$ and
$$F_j[\phi]:=\sigma_k^{1/k}[\lam(-\D^2\phi
   +\frac{n}{n-2}\frac{\D v_j\times \D v_j}{v_j}
   -\frac {1}{n-2}\frac{|\D v_j|^2}{v_j}g_0
   +\frac {n-2}{2} v_jA_{g_0}]
   \ge v_j^p $$
when $\eps>0$ is sufficiently small, where $V_{v_0}$ is the matrix
relative to $v_0$, given in (2.2).  Hence by the concavity of
$\sigma_k^{1/k}$,
$$F^{ab}[v_j]D_{ab}(\phi-v_j)\le F(v_j)-v_j^p\to 0\tag 4.15$$
in $L^{\wtt p}(\Om)$ for any $\wtt p<\infty$, where
$F^{ab}[v_j]=\frac{\p}{\p r_{ab}}\sigma^{1/k}_k(\lam(r))$ at
$r=V_{v_j}$ ($a, b=1, \cdots, n$), which satisfy $\det F^{ab}\ge
C>0$ for some $C>0$ depending only on $n, k$. Applying the
Aleksandrov-Bakelman maximum principle [GT] to (4.15) in
$\{\phi<v_j\}$ and sending $j\to\infty$, we conclude that $\phi\ge
v_0$ near $0$, which is a contradiction so that $F_0[v_0]\le
v_0^p$ at $x_0$. By a similar argument we obtain the reverse
inequality and hence we conclude (4.14) for $v_0$. Since the limit
equation (4.14) is locally uniformly elliptic with respect to
$v_0$, we then conclude further regularity by the Evans-Krylov
estimates and linear theory [GT]. In particular we obtain $v_0\in
C^\infty$. $\square$

A more general approach to the approximation argument to obtain
(4.14) can be obtained by extending the theory of Hessian measures
in [TW1, TW3] to the operators of the type $F_0$.

\proclaim{Lemma 4.3} There exists a function $v_0\in\Phi_k$ such
that the solution $v$ of (4.1) satisfies $J(v(\cdot, t))\ge -C$
and $\sup v(\cdot, t)\ge c_0>0$ for all $t\ge 0$.
\endproclaim

\noo{\it Proof}.  Let $P$ be the set of paths introduced in \S2.2.
For $\gamma\in P$, let $v_s$ ($s\in [0, 1]$) be the solution of
(4.1) with initial condition $v_s(\cdot, 0)=\gamma(s)$. Then by
(4.5) and the comparison principle, there is an $s_0>0$ such that
$v_s(\cdot, t)\to 0$ uniformly for $s\le s_0$. Denote by
$I_\gamma$ the set of $s\in [0, 1]$ such that $J(v_s(\cdot, t))\ge
0$ for all $t>0$. Then $(0, s_0)\subset I_\gamma$.  Let $s^*=\sup
\{s\ |\ s\in I_\gamma\}$.

Obviously $s^*\in I_\gamma$. For if there exists $t$ such that
$J(v_{s^*}(\cdot, t))< 0$, then $J(v_{s}(\cdot, t))< 0$ for
$s<s^*$ sufficiently close to $s^*$, which implies $s^*\ne \sup
\{s\ |\ s\in I_\gamma\}$.  It is also easy to see that
$v_{s^*}(\cdot, t)$ does not converges to zero uniformly, for
otherwise $v_s(\cdot, t)\to 0$ uniformly for $s>s^*$ and near
$s^*$. Finally by our definition of the set $P$, we have $1\not\in
I_\gamma$, namely $s^*<1$. Hence $v_0=\gamma(s^*)$ satisfies Lemma
4.3. $\square$

\vskip10pt

>From the above three Lemmas, one sees that there is a sequence
$t_j\to\infty$ such that $v_{s^*}(\cdot, t_j)$ converges to a
solution of (2.5) for $1<p< \frac{n+2}{n-2}$. Next we prove

\proclaim{Lemma 4.4} For any given $1<p<\frac {n+2}{n-2}$, the set
of solutions of (2.5) is compact.
\endproclaim

\noo{\it Proof}. By the a priori estimates it suffices to show
that the set of solutions is uniformly bounded. If on the contrary
that there is a sequence of solutions $v_j\in\Phi_k$ such that
$\sup v_j\to\infty$, denote $m_j=\sup v_j$ and assume that the sup
is attained at $z_j$. Similar to (4.7) we make a translation and a
dilation of coordinates and a scaling for solution, namely
$$\align
\wtt v_j(y) & =m_j^{-1}v_j(x),\\
y & = R_j x\ \ \ R_j=m_j^{\frac{2}{n-2}-\frac{\eps}{2k}}. \\
\endalign $$
Then $0<\wtt v_j\le 1$, and $\wtt v_j$ satisfies
$$\sigma_k(\lam(\wtt V))=\wtt v^{k\frac {n+2}{n-2}-\eps} .$$
By the a priori estimates in \S 3.1, $\wtt v$ is locally uniformly
bounded in $C^3$. Hence $\wtt v_j$ converges by a subsequence to a
positive solution $\wtt v$ of
$$\sigma_k(\lam(V))=v^{k\frac {n+2}{n-2}-\eps}
                          \ \ \text{in}\ \ \R^n. \tag 4.16$$
By the Liouville Theorem [LL2], there is no nonzero solution to
the above equation. We reach a contradiction. $\square$

\vskip5pt

Let $v$ be a $k$-admissible solution of (2.5). Then we have
$$\int_{(\M, g_0)} v^{\frac {2n-4k}{n-2}}\sigma_k(\lam(\frac Vv))
 -\int_{(\M, g_0)} v^{\frac {2n-4k}{n-2}}v^{\frac {4k}{n-2}-\eps}=0.$$
Hence
$$J(v)=\sup_{t>0} J(tv). $$
By (4.3) we have
$$\align
J(v) & =(\frac {n-2}{2n-4k}-\frac {1}{p+1})
  \int_{(\M, g_0)} v^{\frac {2n-4k}{n-2}}v^{\frac {4k}{n-2}-\eps}\\
  &\ge C>0\tag 4.17\\
  \endalign   $$
By the compactness in Lemma 4.4, the constant $C$ is bounded away
from zero.

\proclaim{Lemma 4.5} There exists a solution $v_p$ of (2.5) such
that $J(v_p)=c_p$.
\endproclaim

\noo{\it Proof}. For any given constant $\delta>0$, choose a path
$\gamma\in P$ such that $\sup_{s\in (0, 1)}J(\gamma(s))\le
c_p+\delta$. By the proof of Lemma 4.3, there exists $s^*\in (0,
1)$ such that the solution of (4.1) with initial condition
$v(\cdot, t)=\gamma_{s^*}$ converges to a solution $v^*_\delta $
of (2.5). Since (4.1) is a descent gradient flow, we have
$J(v^*_\delta)<c_p+\delta$. Letting $\delta\to 0$, by the
compactness in Lemma 4.4, $v^*_\delta$ converges along a
subsequence to a solution $v$ of (2.5) with $J(v)\le c_p$. Note
that $J(v)=\sup_{s>0} J(sv)\ge c_p$. Hence $J(v)=c_p$. $\square$

>From (4.17) we also have
$$c_p\ge C>0. \tag 4.18$$
We have thus proved Theorem 4.1.

\vskip10pt

{\bf 4.2. Proof of Lemma 2.1} (for $2\le k<\frac n2$). Let $v_p$
be a solution of (2.5) with $J_p(v_p)=c_p$. If there is a sequence
$p_j\nearrow \frac {n+2}{n-2}$ such that $\sup v_{p_j}$ is
uniformly bounded, by the a priori estimate in \S 3.1, $v_{p_j}$
sub-converges to a solution of (2.1) and Lemma 2.1 is proved.

If $\sup v_p\to\infty$ as $p\nearrow \frac {n+2}{n-2}$, noting
that $c_p\le \sup_{s>0} J(sv_0)$ for any given $v_0\in\Phi_k$, we
see that $c_p$ is uniformly bounded from above for $p\in [1,
\frac{n+2}{n-2}]$. By (4.17),
$$\int_{(\M, g_0)}v_p^{p+1}\le C,  \tag 4.19$$
where $C$ is independent of $p\le \frac{n+2}{n-2}$. Denote
$m_p=\sup v_p$ and assume that the sup is attained at $z_p=0$. As
before we make a dilation of coordinates and a scaling for
solution, namely
$$\align
\wtt v_p(y) & =m_p^{-1}v_p(x), \\
y & = R_p x\ \ \ R_p=m_p^{\frac{2}{n-2}-\frac{\eps}{2k}}. \\
\endalign $$
Then $0<\wtt v_p\le 1$, and $\wtt v_p$ satisfies
$$\sigma_k(\lam(\wtt V))=\wtt v^{k\frac {n+2}{n-2}-\eps}$$
in $B_{cR_p}$ for some constant $c>0$ independent of $p$. Note
that in the present case, $\eps=\frac{n+2}{n-2}-p\to 0$. By the a
priori estimates in \S 3.1, $\wtt v$ is locally uniformly bounded
in $C^3$. Hence $\wtt v_p$ converges by a subsequence to a
positive solution $\wtt v$ of
$$\sigma_k(\lam(V)=v^{k\frac {n+2}{n-2}}
                          \ \ \text{in}\ \ \R^n. $$
By the Liouville Theorem [LL1],
$$\wtt v(y)=\ol c (1+|y|^2)^{\frac {2-n}{2}}, \tag 4.20$$
where $\ol c=[n(n-2)]^{(n-2)/4}$. Moreover
$$\sup_{s>0} J_{p^*}(s\wtt v; \R^n)=c_{p^*}[S^n] ,\tag 4.21$$
with $p^*=\frac {n+2}{n-2}$, where $c_p$ was defined in (2.9).

The above argument implies that the metric
$g=v_p^{\frac{4}{n-2}}g_0$ is a bubble near the maximum point
$z_p$. By (4.20), $v_p$ has the asymptotical behavior
$$v_p(x)=\ol c \, (\frac {\delta}{\delta^2+r^2})^{\frac {n-2}{2}} (1+o(1))
 \ \ \ \ \delta=m_p^{-\frac 2{n-2}+\frac{\eps}{2k}}. \tag 4.22$$
For a sufficiently small $\th>0$, let $\Om_p=\{x\in\M\ |\
v_p(x)>\th m_p\}$, and let
$$\hat v_p(x)=\cases
v_p(x)\ \ \ x\in\M-\Om_p,\\
\th m_p\ \ \ x\in\Om_p.\\
\endcases$$
Note that by assumption (1.7),
$$\sup_{s>0}J(s v_p)=c_p<c_{p^*}[S^n]\tag 4.23$$
when $p<\frac{n+2}{n-2}$ and is close to $\frac {n+2}{n-2}$.

Combining (4.21), (4.22) and (4.23), we see that
$$\int_{(\M, g_0)} {\hat v}_p^{p+1}\ge C>0$$
for some $C$ independent of $\th$, provided $\th$ is sufficiently
small and $m_p$ is sufficiently large, and
$$\sup_{s>0} J(s\hat v_p)<\sup_{s>0}J(s v_p).$$
Namely $\sup_{s>0} J(s\hat v_p)<c_p$, which is in contradiction of
our definition of $c_p$. Note that $\hat v_p$ is not smooth, but
can be approximated by smooth, $k$-admissible functions.  This
completes the proof of Lemma 2.1. $\square$

\vskip10pt

{\bf 4.3. The case $k=\frac n2$}. In this case, the proof of Lemma
4.3 does not apply,  due to that $J(v)\to -\infty$ as $v\to 0$,
and also we don't know if $\E_{n/2}(v)$ is bounded from below for
any admissible function $v$ with Vol$\M_{g_v}=1$, where
$$J(v)=\E_{n/2}(v)
 -\frac {1}{p+1}\int_{(\M, g_0)} v^{\frac{2n}{n-2}-\eps}$$
is the corresponding functional and $\E_{n/2}$ is given in (2.32).
However when $k=\frac n2$, we have the following

\proclaim{Lemma 4.6} Assume that equation (2.1) is variational.
Then $\F_{n/2}(v)$ is a constant.
\endproclaim

\noo{\it Proof}. When $k=\frac n2$, we write the equation (2.1)
and the functional $\F_{n/2}$ in the form
$$\align
\sigma_{n/2}(\lam(W))& =e^{-nw}, \\
\F_{n/2}(w)& =\int_{(\M, g_0)} \sigma_{n/2}(\lam(W)),\tag 4.24\\
\endalign $$
To prove that $\F_{n/2}$ is equal to a constant, we have
$$\align
\F_{n/2}(w)-\F_{n/2}(w_0)
 &=\int_{(\M, g_0)} \int_0^1\frac{d}{dt}\sigma_{n/2}(\lam(W_t))\\
 &=\int_0^1\int_{(\M, g_0)} L_{w_t}(w)\\
 \endalign $$
where $w_t=tw$, $w_0=0$, $L_{w_t}$ is the linearized operator of
$\sigma_{n/2}(\lam(W))$ at $w_t$. By the assumption that equation
(2.1) is variational, we have (see \S2.4)
$$\int_{(\M, g_0)} L_{w_t}(w)=\int_{(\M, g_0)} w L_{w_t}(1)=0. $$
This completes the proof of Lemma 4.6. $\square$

By assumption (1.7), we have
$$\F_{n/2}(v)=c_0<Y_{n/2}(S^n)\tag 4.25$$
for some constant $c_0$ depending on $(\M, g_0)$. Lemma 4.6
enables us to prove the following

\proclaim{Lemma 4.7} For $1<p\le \frac {n+2}{n-2}$, the set of
solutions of (2.5) is compact.
\endproclaim

\noo{\it Proof}. When $1<p<\frac {n+2}{n-2}$, the proof is the
same as that of Lemma 4.4.

When $p=\frac {n+2}{n-2}$, we use the same argument of Lemma 4.4.
Instead of (4.16), we have the equation
$$\sigma_{n/2}(\lam(V))=v^{\frac n2\frac {n+2}{n-2}}
                         \ \ \text{in}\ \ \R^n. \tag 4.26$$
By the Liouville theorem [LL1], $v$ must be the function given in
(4.20). Hence we have
$$\int_{\R^n} \sigma_{n/2}(\lam(V))= Y_{n/2}(S^n).$$
By (4.9), we obtain that
$$\underline\lim_{j\to\infty}\F_{n/2}(v_j)\ge Y_{n/2}(S^n). $$
This is in contradiction with (4.25). $\square$

By Lemma 4.7, we can prove the existence of solutions of (2.1) by
a degree argument, see [CGY2, LL1]. We omit the details here.

\vskip10pt

{\bf 4.4. A Sobolev type inequality}. As a consequence of our
argument above, we have the following Sobolev type inequality.

\proclaim{Corollary 4.1} Let $2\le k<\frac n2$. Then there exists
a constant $C>0$ such that the inequality
$$\big[ Vol(\M_g)]^{\frac {n-2}{2n}}\le C\big[\int_\M
\sigma_k(\lam(A_g))d\, vol_g]^{\frac {n-2}{2n-4k}} \tag 4.27$$
holds for any conformal metric $g=v^{\frac {4}{n-2}}g_0$ with
$v\in\Phi_k$.
\endproclaim

\noo{\it Proof}. Note that (4.27) is equivalent to
$$\big[\int_{(\M, g_0)} v^{\frac{2n}{n-2}} \big]^{\frac{n-2}{2n}}
 \le C \big[ \int_{(\M, g_0)} v^{\frac{2n-4k}{n-2}}
      \sigma_k(\lam(\frac V v))\big]^{\frac {n-2}{2n-4k}},\tag 4.28$$
which is equivalent to (4.18). $\square$

\noo{\bf Remark 4.1}. The Sobolev type inequality (4.27) is
similar to (2.17) of $l=0$, $k\ge 1$, and was proved in [GW2] for
locally conformally flat manifolds. If $\M$ is locally conformally
flat and $f=\mu(e^{-2kw})$, the flow (3.21) has a remarkable
property. That is by the moving plane argument of Ye [Ye] and the
conformal invariance of the equation, one obtains the gradient
estimate (3.23) at all time $t$, with the upper bound $C$
depending only on the initial function $w(\cdot, 0)$, for any
monotone increasing $\mu$ satisfying (2.13). Therefore by the
second derivative estimate (3.26), the solution of (3.21)
converges to a solution of the $\kk$-Yamabe equation (2.1), and
accordingly one also obtains the Sobolav type inequality (4.27).
Theorem 4.2 shows that the Sobolev type inequality also holds on
general manifolds provided equation (1.1) is variational.

\vskip30pt

\newpage

\centerline{\bf 5. Proof of Lemma 2.2}

\vskip10pt

We let $v_\eps$ be the function given by
$$v_\eps(x)=
  \bigg(\frac{\eps}{\eps^2+r^2}\bigg)^{\frac{n-2}{2}},\tag 5.1$$
where $r=|x|$, $x\in\R^n$, and $\eps>0$ is a small constant. Let
$V_\eps$ be the matrix relative to $v_\eps$, see (2.2). Then we
have
$$\frac {V_\eps}{v_\eps}=(n-2)v_\eps^{\frac {4}{n-2}}I. $$
Hence $v_\eps$ is $k$-admissible on $\R^n$ and
$$\sigma_k(\lam(\frac {V_\eps}{v_\eps})
                          =C_{n, k} v_\eps^{\frac{4k}{n-2}}, $$
where $C_{n, k}=\frac{n!(n-2)^k}{k!(n-k)!}$. It follows that
$$\int_{\R^n} v^{\frac {2n-4k}{n-2}}\sigma_k(\lam(\frac Vv))
 =C_{n, k}\int_{R^n} v^{\frac{2n}{n-2}}. $$
So we have
$$Y_k(S^n)= \frac{\int_{\R^n} v^{\frac
{2n-4k}{n-2}}\sigma_k(\lam(\frac Vv))}
 {\big[\int_{\R^n} v^{\frac {2n}{n-2}}\big]^{(n-2k)/n}}
  =C_{n, k}\bigg[\int_{\R^n} v^{\frac{2n}{n-2}}\bigg]^{2k/n}.\tag 5.2$$
In particular we have
$$Y_k(S^n)=\frac{C_{n, k}}{(n(n-2))^k} [Y_1(S^n)]^k. \tag 5.3$$
In the above $v=v_\eps$ and the integrations are independent of
$\eps$.

To verify (1.7), it would be natural to use the function (5.1) as
a test function, as in the case $k=1$ [Au1, S1]. However on a
general manifold, the function $v_\eps$, where $r$ denotes the
geodesic distance from a given point, is $k$-admissible only when
$r\le C\eps^{1/2}$. It seems impossible to find an explicit test
function.

Instead we shall deduce (1.7) directly from (1.6). First note that
by assumption, there exists a function $v>0$ such that $\wtt
g=v^{4/(n-2)}g\in\Ga_k$. Hence $\sigma_1(\lam(A_{\wtt g}))>0$.
That is the scalar curvature of $(\M, \wtt g)$ is positive. Hence
the comparison principle for the operator $\sigma_1(\lam(A_g))$
holds on $\M$.

Let $v_1$ be a solution of the Yamabe problem (with $k=1$) such
that $Q_1(v_1)<Y_1(S^n)$, where
$$Q_1(v)=\frac{\int_\M v\sigma_1(\lam(V))}
              {[\int_\M v^{2n/(n-2)}]^{(n-2)/n}}. $$
Let $v_k$ be the solution of
$$\sigma_k(\lam(V))= C_{n, k} v_1^{k\frac {n+2}{n-2}}
                           \ \ \ \text{in}\ \ \M. \tag 5.4$$
By Lemma 2.3(vi),  we have
$$ -\Delta v_k+\frac{n-2}{4(n-1)}Rv_k
     =\sigma_1(\lam(V_k)) \ge  n(n-2) v_1^{\frac{n+2}{n-2}}. $$
Since $v_1$ satisfies
$$ -\Delta v+\frac{n-2}{4(n-1)}Rv
    =  n(n-2) v_1^{\frac{n+2}{n-2}}, $$
by the comparison principle,
$$v_k\ge v_1\tag 5.5$$
Now writing
$$Q_k(v)=\frac
 {\int_\M v^{\frac {2n}{n-2}-k\frac{n+2}{n-2}} \sigma_k(\lam(V))}
 {\big[ \int_\M v^{2n/(n-2)}\,d  vol_{g}\big]^{(n-2k)/n}},$$
we claim that
$$Q_k(v_k) < Y_k(S^n),\tag 5.6$$
namely (1.7) holds. Indeed, when $k\ge 2$, we have $\frac
{2n}{n-2}-k\frac{n+2}{n-2}<0$. Hence by (5.5),
$$v_1^{\frac {2n}{n-2}-k\frac{n+2}{n-2}}
          \ge v_k^{\frac{2n}{n-2}-k\frac{n+2}{n-2}}.$$
Hence
$$\align
\int_\M v_k^{\frac {2n}{n-2}-k\frac{n+2}{n-2}}
                           \sigma_k(\lam(V_k)) \, d vol_{g}
&\le C_{n, k}\int_{B_\rho} v_1^{\frac {2n}{n-2}-k\frac{n+2}{n-2}}
                 v_1^{k\frac{n+2}{n-2}} \, d vol_{g}\\
& \le C_{n, k}\int_{B_\rho} v_1^{\frac {2n}{n-2}}\, d vol_{g}\\
\endalign $$
and
$$\int_\M v_k^{\frac {2n}{n-2}}\, d vol_{g}
     \ge \int_{B_\rho} v_1^{\frac {2n}{n-2}}\, d vol_{g}. $$
Therefore we obtain
$$Q_k(v_k)
 \le C_{n, k} \big[\int_\M v_1^{\frac {2n}{n-2}}\,
                 d vol_{g}\big]^{2k/n},\tag 5.7$$
so that (5.6) follows from (5.3).

We remark that a similar argument can be used to prove the
inequality $Y_1(\M)<Y_1(S^n)$ for some manifolds.

\vskip30pt

\newpage

\centerline{\bf 6. Some remarks}

\vskip10pt

{\bf 6.1. Compactness of the solution set}. For the Yamabe problem
$(k=1)$, Schoen [S2] has shown furthermore that the set of
solutions is compact if the manifold is locally conformally flat
and not conformally equivalent to the sphere. Schoen's result was
extended to general compact manifolds for which the positive mass
theorem holds, such as in low dimensions $3\le n\le 7$ [LZ, M].

When $k>\frac n2$, the compactness of solutions has also been
established in [GV2], for the more general equation
$$\sigma_k(\lam(A_{g_v}))=f\, v^{\frac{n+2}{n-2}},\tag 6.1$$
where $f$ is any positive, smooth function $f$. Their proof relies
crucially on the fact that the Ricci curvature $Ric_{g_v}\ge C
g_v$ if $g_v=v^{\frac{4}{n-2}}g_0$ is a solution to (6.1). For the
convenience of the reader we indicate the idea of their proof
here. From the positivity of the Ricci curvature, the volume of
$(\M, g_v)$ is uniformly bounded. Hence if there exists a sequence
of solutions $\{v_k\}$ with $\sup v_k\to\infty$, there are at most
finitely many blow-up points $P=:\{p_0, p_1, \cdots, p_s\}$. By
the interior first and second derivative estimates (3.23) and
(3.26), $v_k/\inf v_k$ converges locally uniformly on $\M-P$ to a
$C^{1,1}$ function $v$ with $\sigma_k(\lam(A_{g_v}))=0$. By a
technical analysis one has $v(x)=C(1+o(1))r^{2-n}$ near the
singularity set $P$. Since $Ric_{g_v}>0$, the ratio $h(r)=:\frac
{|B_r(p_0)|}{r^n}$ is non-increasing, where $|B_r(p_0)|$ denotes
the volume of the geodesic ball on $(\M, g_v)$. On the other hand,
since $v(x)=C(1+o(1))r^{2-n}$, we have $\lim_{r\to\infty}
h(r)=(s+1)\omega_n$, where $\omega_n$ is the volume of the
Euclidean unit ball. Hence $s=0$ and $(\M, g_v)$ is isometric to
the Euclidean space $\R^n$, and so $\M$ is conformal to the unit
sphere. Note that when $2\le k\le \frac n2$, the Ricci curvature
of $(\M, g_v)$ may not be positive anymore.

\vskip10pt

{\bf 6.2. Conditions (C1) and (C2)}. As indicated earlier, we
impose condition (C1) so that equation (1.1) is elliptic. If a
fully nonlinear partial differential equation is not elliptic,
little is known about the existence and regularity of solutions.
For example it is unknown whether there is a local solution to the
Monge-Ampere equation $\det D^2 u=f$ when the right hand side $f$
changes sign, even in dimension two. But possibly condition (C1)
may be replaced by the positivity of the Yamabe constant
$Y_k(\M)$, as in the case $k=2$, $n=4$ [CGY1, GV1].

As for the condition (C2), the variational approach to the
$\kk$-Yamabe problem is natural, as in the case $k=1$. Indeed this
approach has already been employed in [CGY1, CGY2, GW2, GW3], and
Theorem 2.1 was proved in [CGY1, CGY2] when $n=4$, $k=2$, and in
[GW2, GW3] when $\M$ is locally conformally flat. At the moment we
are not aware of any other possible ways to remove the variational
structure condition (C2) for the case $2\le k\le \frac n2$, even
in low dimensional cases.

\vskip10pt

{\bf 6.3. The full $k$-Yamabe problem [K, La]}. We bring to the
attention of the readers the full $k$-Yamabe problem. On a
Riemannian manifold $(M^n, g)$, one can define a series of scalar
curvatures
$$s_k=s_k(Riem)=s_k(W +A\odot g),\tag 6.2$$
for $k=1, 2, \cdots, [\frac n2]$, where $Riem, W, A$ are
introduced at the beginning in the introduction. The $k$-scalar
curvature can also be expressed simply as
$$s_k=\frac {1}{(2k)!} c^{2k} Riem^k,\tag 6.3$$
where $c$ is the standard contraction operator, and
$Riem^k=Riem\circ\cdots\circ Riem$ is the product introduced in
[K].

When $k=1$, $s_1$ is the usual scalar curvature. When $n$ is even,
$s_{n/2}$ is the Lipschitz-Killing curvature. Furthermore if $\M$
is a hypersurface, the $k$-scalar curvature $s_k$ is the $2k^{th}$
mean curvature $H_{2k}$, which is equal to the $2k^{th}$
elementary symmetric polynomial of the principal curvatures of the
hypersurface, which is intrinsic quantity [Sp].  When $\M$ is
locally conformally flat, then the Weyl curvature in (6.2)
vanishes, and $s_k$ turns out to be the $k$-curvature given in
(1.1).

The full $k$-Yamabe problem concerns the existence of a conformal
metric such that the $k$-scalar curvature is a constant. This
problem coincides with the $\kk$-Yamabe problem for locally
conformally flat manifolds. The corresponding equation of the
$k$-Yamabe problem is always variational, as in the case $k=1$
[La]. However the equation is of mixed type in general.

\vskip30pt

\newpage

\baselineskip=12.0pt
\parskip=0pt

\Refs\widestnumber\key{ABC}

\item {[Au1]} T. Aubin, Equations diff\'erentielles non lin\'eaires
       et probl\`me de Yamabe concernant la courbure scalaire.
       J. Math. Pures Appl. (9) 55 (1976), 269--296.

\item {[Au2]} T. Aubin,
       Some nonlinear problems in Riemannian geometry,
       Springer, 1998.

\item {[B]} A.L. Besse,
       Einstein manifolds, Springer-Verlag, Berlin, 1987.

\item {[Br]} T. Branson,
       Differential operators canonically associated to a
       conformal structure,
       Math. Scand.,  57(1985), 293--345.

\item {[BV]} S. Brendle and J. Viaclovsky,
       A variational characterization for $\sigma_{n/2}$,
       Calc. Var. PDE. 20(2004), 399--402.

\item {[CNS]} L.A. Caffarelli, L. Nirenberg, and J. Spruck,
       Dirichlet problem for nonlinear second order
       elliptic equations III.
       Functions of the eigenvalues of the Hessian,
       Acta Math. {\bf 1985}, {\it 155}, 261--301.

\item {[CGY1]} A. Chang, M. Gursky, P. Yang,
       An equation of Monge-Am\`ere type in conformal geometry,
       and four-manifolds of positive Ricci curvature,
       Ann. of Math. (2) 155(2002), 709--787.

\item {[CGY2]} A. Chang, M. Gursky, P. Yang,
       An a priori estimate for a fully nonlinear equation on
       four-manifolds,
       J. Anal. Math. 87 (2002), 151--186.

\item {[CGY3]} A. Chang, M. Gursky, and P. Yang,
       Entire solutions of a fully nonlinear equation,
       Lectures on partial differential equations, 43--60,
       Inter. Press, 2003.

\item {[Ch1]} K.S. Chou (K. Tso),
      On a real Monge-Ampere functional,
      Invent. Math. 101(1990), 425--448.

\item {[CW1]} K.S. Chou and X-J. Wang,
       Variational solutions to Hessian equations,
       preprint, 1996.
       Available at wwwmaths.anu.edu.au/research.reports/96mrr.html.

\item {[CW2]} K.S. Chou and X-J. Wang,
       A variational theory of the Hessian equation,
       Comm. Pure Appl. Math. 54 (2001), 1029--1064.

\item {[GV1]} M. Gursky and J. Viaclovsky,
       A fully nonlinear equation on four-manifolds with positive scalar curvature,
       J. Differential Geom. 63 (2003), 131--154.

\item {[GV2]} M. Gursky and J.Viaclovsky,
       Prescribing symmetric functions of the eigenvalues of the
       Ricci tensor, arXiv:math.DG/0409187.

\item {[GW1]} P. Guan and G. Wang,
       Local estimates for a class of fully nonlinear equations
       arising from conformal geometry,
       Int. Math. Res. Not. (2003), 1413--1432.

\item {[GW2]} P. Guan and G. Wang,
        A fully nonlinear conformal flow on locally conformally
        flat manifolds,
        J. Reine Angew. Math. 557 (2003), 219--238.


\item {[GW3]} P. Guan and G. Wang,
       Geometric inequalities on locally conformally flat manifolds,
       Duke Math. J. 124 (2004), 177--212.

\item {[GW4]} P. Guan and Xu-Jia Wang,
       On a Monge-Amp\`ere equation arising in geometric optics,
       J. Diff. Geom., 48(1998), 205--223.

\item {[I1]} N. Ivochkina,
       Solution of the Dirichlet problem for certain equations of
       Monge-Amp\`ere type (Russian), Mat. Sb. (N.S.) 128(170) (1985),
       403--415.

\item {[K]} R.S. Kulkarni,
       On the Bianchi identity,
       Math. Ann., 199(1972), 175-204.

\item {[La]} M.-L. Labbi,
       On a variational formula for the H. Weyl curvature
       invariants, arXiv:math.DG/0406548.

\item {[Lb]} D. Labutin,
       Potential estimates for a class of fully nonlinear elliptic
       equations. Duke Math. J. 111 (2002), 1--49.

\item {[LP]} J.M. Lee and T.H. Parker,
        The Yamabe problem, Bull. AMS 17(1987), 37--91.

\item {[LL1]} A. Li and Y.Y. Li,
       On some conformally invariant fully nonlinear equations.
       Comm. Pure Appl. Math. 56 (2003), 1416--1464.

\item {[LL2]} A. Li and Y.Y. Li,
       On some conformally invariant fully nonlinear equations
       I\!I, Liouville, Harnack, and Yamabe, preprint.

\item {[LS]} J. Li and W.M. Sheng,
      Deforming metrics with negative curvature by a fully nonlinear
      flow, Calc. Var. PDE, to appear.

\item {[LZ]} Y.Y. Li and L. Zhang,
      Compactness of solutions to the Yamabe problem,
      C. R. Math. Acad. Sci. Paris 338 (2004), 693--695.

\item {[LT]} M. Lin and N.S. Trudinger,
      On some inequalities for elementary symmetric functions
      Bull. Aust. Math. Soc., 50(1994), 317--326.

\item {[MTW]} X.N. Ma, N.S. Trudinger, and X-J. Wang,
      Regularity of potential functions to the optimal
      transportation problem, Arch. Rat. Mech. Anal., to appear.

\item {[M]} F.C. Marques,
      A priori estimates for the Yamabe problem in non-locally
      conformal flat cases, \newline arXiv:math.DG/0408063.

\item {[O]} P.J. Olver,
       Applications of Lie groups to differential equations,
       Springer, 1993.

\item {[S1]} R. Schoen,
       Conformal deformation of a Riemannian metric to constant
       scalar curvature,
       J. Diff. Geom. 20(1984), 479--495.

\item {[S2]} R. Schoen,
       Variational theory for the total scalar curvature
       functional for Riemannian metrics and related topics.
       Topics in Calculus of Variations.
       Lectures Notes in Math. 1365. Springer, Berlin (1989), pp.
       120--154.

\item {[SY1]} R. Schoen and S.T. Yau,
       Conformally flat manifolds, Kleinian groups and scalar
       curvature,  Invent. Math. 92 (1988), 47--71.

\item {[SY2]} R. Schoen and S.T. Yau,
       Lectures on Differential geometry.
       International Press, 1994.

\item {[Sc]} F. Schulz,
      Regularity theory for quasi-linear elliptic systems and
      Monge-Amp\`ere equations in two dimensions,
      Lecture Notes in Math. 1445 (1990).

\item {[Sp]} M. Spivak,
       A comprehensive introduction to differential geometry,
       Vol.4, Publish or Perish Inc. 1979.

\item {[Tr1]} N.S. Trudinger,
       Remarks concerning the conformal deformation of Riemannian
       structures on compact manifolds,
       Ann. Scuola Norm. Sup. Pisa (3) 22(1968), 265--274.

\item {[TW1]} N.S. Trudinger and X-J. Wang,
        Hessian measures I,
        Topol. Methods Nonlinear Anal. 10 (1997), 225--239.

\item {[TW2]} N.S. Trudinger and X-J. Wang,
       Hessian measures I\!I,
       Ann. of Math. (2) 150 (1999), 579--604.

\item {[TW3]} N.S. Trudinger and X-J. Wang,
       Hessian measures I\!I\!I,
       J. Funct. Anal. 193 (2002), 1--23.

\item {[TW4]} N.S. Trudinger and X-J. Wang,
       A Poincar\'e type inequality for Hessian integrals,
       Calc. Var. PDE  6(1998), 315--328.

\item {[TW5]} N.S. Trudinger and X-J. Wang,
       On the weak continuity of elliptic operators and applications
       to potential theory, Amer. J. Math. 124 (2002), 369--410.

\item {[V1]} J. Viaclovsky, Conformal geometry, contact geometry,
      and the calculus of variations.
      Duke Math. J. 101 (2000), no. 2, 283--316.

\item {[V2]} J. Viaclovsky,
       Estimates and existence results for some fully nonlinear
       elliptic equations on Riemannian manifolds,
       Comm. Anal. Geom. 10 (2002), 815--846.

\item {[W1]} X-J. Wang,
      A class of fully nonlinear elliptic equations and related
      functionals, Indiana Univ. Math. J. 43 (1994), 25--54.

\item {[W2]} X-J. Wang,
      On the design of a reflector antenna,
      Inverse Problems, 12(1996), 351--375.

\item {[Ya]} H. Yamabe,
      On a deformation of Riemannian structures on compact manifolds,
      Osaka Math. J. 12(1960), 21--37.

\item {[Ye]} R. Ye,
       Global existence and convergence of Yamabe flow.
       J. Differential Geom. 39 (1994),  35--50.

\endRefs

\enddocument
\end